\documentclass[review]{elsarticle}
\usepackage[OT1]{fontenc}  	
\usepackage[utf8]{inputenc}
\usepackage[spanish,french,english,USenglish]{babel} 		
\usepackage{amssymb}  			
\usepackage{amsmath}	
\usepackage{amsthm}	
\usepackage{esint}  				
\usepackage{bbm}   				
\usepackage{color,graphicx} 			
\usepackage[small]{caption}
\usepackage{cases}	
\usepackage{bigints_mod}	
\usepackage{mathpazo}

\newcommand{\be}{\begin{equation}}
\newcommand{\ee}{\end{equation}}

\renewcommand{\Re}{\mathop{\rm Re}\nolimits}

\newcommand{\tg}{\mathop{\rm tg}\nolimits}

\newcommand{\arctg}{\mathop{\rm arctg}\nolimits}

\newcommand{\res}{\mathop{\rm res}\limits}

\makeatletter
\newcommand{\specialnumber}[1]{
\def\tagform@##1{\maketag@@@{(\ignorespaces##1\unskip\@@italiccorr#1)}}}
\makeatother

\makeatletter
\def\ps@pprintTitle{%
\let\@oddhead\@empty
\let\@evenhead\@empty
\let\@oddfoot\@empty
\let\@evenfoot\@oddfoot
}\makeatother

\journal{Journal of Number Theory (Elsevier)}

\begin{document}

\begin{frontmatter}

\title{Expansions of generalized Euler's constants into the series of
polynomials in $\pi^{-2}$
and into the formal enveloping series \\ with rational coefficients only\footnote{\phantom{-} \\[3mm]
\texttt{\underline{Note to the readers of the 4th arXiv version:}
this version is a copy of the journal version of the article, which has been published in the Journal of Number Theory (Elsevier), vol.~158, pp.~365-396, 2016.
DOI 10.1016/J.JNT.2015.06.012  http://www.sciencedirect.com/science/article/pii/S0022314X15002255 \\
Artcile history: submitted 1 January 2015, accepted 29 June 2015, published on-line 18 August 2015.\\ 
The layout of the present version and its page numbering differ from the journal version, but the content, the numbering of equations and the numbering 
of references are the same. This version also incorporates some minor corrections to the final journal version, which were published on-line
in the same journal on December 8, 2016 (DOI 10.1016/J.JNT.2016.11.002).
For any further reference to the material published here, please, use the journal version of the paper, 
which you can always get for free by writing a kind e-mail to the author.\\[-16mm]}}}

\author{Iaroslav V.~Blagouchine\corref{cor1}} 
\ead{iaroslav.blagouchine@univ-tln.fr, iaroslav.blagouchine@centrale-marseille.fr}
% \ead[url]{http://lseet.univ-tln.fr/\~{}iaroslav/}
\cortext[cor1]{Corresponding author. Phones: +33--970--46--28--33, +7--953--358--87--23.}
\address{University of Toulon, France.}

\begin{abstract}
In this work, two new series expansions for generalized Euler's
constants (Stieltjes constants) $\gamma_m$ are obtained.
The first expansion involves Stirling numbers of the first kind,
contains polynomials in $\pi^{-2}$ with rational coefficients and converges
slightly better than Euler's series $\,\sum n^{-2}$.
The second expansion is a semi-convergent series with rational
coefficients only.
This expansion is particularly simple and involves Bernoulli numbers
with a non-linear
combination of generalized harmonic numbers. It also permits to derive
an interesting estimation for generalized Euler's constants, which is
more accurate than several well-known estimations. Finally,
in Appendix A, the reader will also find two simple integral
definitions for the Stirling numbers of the first kind, as well
an upper bound for them.
\end{abstract}

\begin{keyword}
{Generalized Euler's constants},
{Stieltjes constants},
{Stirling numbers},
{Factorial coefficients},
{Series expansion},
{Divergent series},
{Semi-convergent series},
{Formal series},
{Enveloping series},
{Asymptotic expansions},
{Approximations},
{Bernoulli numbers},
{Harmonic numbers},
{Rational coefficients},
{Inverse pi}.
\end{keyword}
\end{frontmatter}
%
%spell_from *************** Text entry area ******************%

%%% PASTABA: lygtys centre

%s1 #&#
\section{Introduction and notations}
%s1.1 #&#
\subsection{Introduction}
The $\zeta$-function, which is usually introduced via one of the
following series,
\be\label{kj023dndndr3}
\zeta(s)=
\begin{cases}
\displaystyle\sum\limits_{n=1}^\infty \frac{1}{\,n^{s}}\,,\qquad & \Re{s}>1 \\[6mm]
\displaystyle\frac{1}{1-2^{1-s}}\sum\limits_{n=1}^\infty \frac{(-1)^{n-1}}{\,n^{s}}\,,\qquad & \Re{s}>0 \,,\quad s\neq1
\end{cases}
\ee
is of fundamental and long-standing importance in modern analysis, number theory, theory of special
functions and in a variety other fields. It is well known that $\zeta(s)$ is meromorphic on the entire complex
$s$-plane and that it has one simple pole at $s=1$ with residue 1.
Its expansion in the Laurent series in a neighbourhood of $s=1$ is
usually written the following form
\be\label{dhd73vj6s1}
\zeta(s)\,=\,\frac{1}{\,s-1\,} + \sum_{m=0}^\infty \frac{(-1)^m (s-1)^m}{m!} \gamma_m\,, 
\qquad\qquad \qquad s\neq1.
\ee
where coefficients $\gamma_m$, appearing in the regular part of
expansion {\eqref{dhd73vj6s1}}, are called
\emph{generalized Euler's constants} or \emph{Stieltjes constants},
both names being in use.\footnote{The definition
of Stieltjes constants accordingly to formula {\eqref{dhd73vj6s1}} is
due to Godfrey H.~Hardy.
Definitions, introduced by Thomas Stieltjes and Charles
Hermite between 1882--1884, did not contain coefficients $(-1)^m$ and $m!$
In fact, use of these factors is not well justified; notwithstanding,
Hardy's form {\eqref{dhd73vj6s1}} is largely accepted and is more
frequently encountered in modern literature.
For more details, see \cite[vol.~I, letter 71 and following]{stieltjes_01},
\cite[p.~562]{lagarias_01}, \cite[pp.~538--539]{iaroslav_07}.}\up{,}\footnote{Some authors use the name \emph{generalized Euler's constants}
for other constants,
which were conceptually introduced and studied by Briggs in 1961 \cite{briggs_01} and Lehmer in 1975 \cite{lehmer_01}. They were subsequently
rediscovered in various (usually slightly different) forms by several
authors, see e.g.~\cite{tasaka_01,pilehrood_01,xia_01}.
Further generalization of both, generalized Euler's constants defined
accordingly to {\eqref{dhd73vj6s1}} and generalized Euler's constants introduced
by Briggs and Lehmer, was done by Dilcher in \cite{dilcher_01}.}
Series {\eqref{dhd73vj6s1}} is the standard definition for $\gamma_m$.
Alternatively, these constants may be also defined via the following limit
\begin{eqnarray}
\label{k98y9g87fcfcf} 
\gamma_m = \lim_{n\to\infty} \left\{
\sum_{k=1}^n \frac{\ln^m k}{k} -
\frac{\ln^{m+1} n}{m+1} \right\} , \quad m=0, 1, 2,\ldots
\end{eqnarray}
The equivalence between definitions {\eqref{dhd73vj6s1}} and {\eqref{k98y9g87fcfcf}} 
was demonstrated by various authors, including Adolf Pilz \cite{gram_01}, Thomas Stieltjes, 
Charles Hermite \cite[vol.~I, letter 71 and following]{stieltjes_01}, Johan Jensen \cite{jensen_02,jensen_03}, 
J\'er\^ome Franel \cite{franel_01},
J{\o}rgen P.~Gram \cite{gram_01},
Godfrey H.~Hardy \cite{hardy_03}, Srinivasa Ramanujan
\cite{ramanujan_01}, William E.~Briggs, S.~Chowla
\cite{briggs_02} and many others, see e.g.~\cite{berndt_02,todd_01,israilov_01,zhang_01}.
It is well known that $\gamma_0=\gamma$ Euler's constant, see
e.g.~\cite{zhang_01}, \cite[Eq.~(14)]{iaroslav_07}.
Higher generalized Euler's constants are not known to be related to the
``standard'' mathematical constants,
nor to the ``classic'' functions of analysis.

In our recent work \cite{iaroslav_08}, we obtained two interesting
series representations for the logarithm of the $\Gamma$-function
containing Stirling numbers of the first kind $S_1(n,k)$
\begin{eqnarray}
\displaystyle\notag
\ln\Gamma(z)\, =&& \left(z-\frac{1}{\,2\,}\right)\!\ln z -z +\frac{1}{\,2\,}\ln2\pi + \\[1mm]
\displaystyle&&\displaystyle\label{kj20ejcn2dnd}
\qquad\qquad
+ \frac{1}{\,\pi\,}\!
\sum_{n=1}^\infty\!\frac{ 1}{\,n\cdot n!\,}\sum_{l=0}^{\lfloor\!\frac{1}{2}n\!\rfloor} (-1)^{l} 
\frac{\, (2l)!\cdot|S_1(n,2l+1)|\,}{(2\pi z)^{2l+1}} \\[5mm]
\displaystyle
\ln\Gamma(z) = &&\displaystyle\left(z-\frac{1}{\,2\,}\right)\!\ln\!{\left(z-\frac{1}{\,2\,}\right)}
- z +\frac{1}{\,2\,}+\frac{1}{\,2\,}\ln2\pi - \notag\\[1mm]
\displaystyle&&\displaystyle\label{lj023od230dend}
\qquad\qquad
- \frac{1}{\,\pi\,}\!
\sum_{n=1}^\infty\!\frac{ 1}{\,n\cdot n!\,}\sum_{l=0}^{\lfloor\!\frac{1}{2}n\!\rfloor} (-1)^{l} 
\frac{\, (2l)!\cdot(2^{2l+1}-1)\cdot|S_1(n,2l+1)|\,}{(4\pi)^{2l+1}\cdot\big(z-\frac{1}{2}\big)^{2l+1}}
\end{eqnarray}
as well as their analogs for the polygamma functions $\Psi_k(z)$.\footnote{Both series converge in a part of the
right half--plane \cite[Fig.~2]{iaroslav_08} at the
same rate as $\sum \big(n\ln^{m}\!
n\big)^{-2}\,$, where $m=1$ for $\ln\Gamma(z)$
and $\Psi(z)$, 
$m=2$ for $\Psi_1(z)$ and $\Psi_2(z)$,
$m=3$ for $\Psi_3(z)$ and $\Psi_4(z)$, \emph{etc.} \label{gtf1a}}
The present paper is a continuation of this previous work, in which we show that
the use of a similar technique permits to derive two new series
expansions for
generalized Euler's constants $\gamma_m$, both series involving
Stirling numbers of the first kind.
The first series is convergent and contains polynomials in $\pi^{-2}$
with rational coefficients (the latter involves Stirling numbers of the
first kind). From this series,
by a formal procedure, we deduce the second expansion, which is
semi-convergent and contains rational terms only.
This expansion is particularly simple and involves only Bernoulli
numbers and a non-linear
combination of generalized harmonic numbers.
Convergence analysis of discovered series shows that the
former converges slightly better than Euler's series
$\sum n^{-2}$, in a rough approximation at the same rate as
\be
\nonumber
\sum_{n=3}^\infty\frac{\ln^m \!\ln n}{\,n^2\ln^2 \! n\,}\,,
\qquad\quad m=0, 1, 2,\ldots
\ee
The latter series diverges very quickly, approximately as
\be
\nonumber
\sum_{n=2}^\infty(-1)^{n-1} \frac{\,\ln^m \!n\,}{\sqrt{n\,}}\left
(\frac{n}{\pi e}\right)^{2n}\,,
\qquad\quad m=0, 1, 2,\ldots
\ee

%s1.2 #&#
\subsection{Notations and some definitions}\label{notations}
Throughout the manuscript, following abbreviated notations are used: $
\gamma=0.5772156649\ldots$ for
Euler's constant, $\gamma_m$ for $m$th generalized Euler's constant
(Stieltjes constant)
accordingly to their definition {\eqref{dhd73vj6s1}},\footnote{In
particular $\gamma_1=-0.07281584548\ldots{}$,
$\gamma_2=-0.009690363192\ldots{}$, $\gamma_3=+0.002053834420\ldots{}$.}
$\binom{k}{n}$ denotes the binomial coefficient $C^n_k$,
${B}_n$~stands for the $n$th Bernoulli number,\footnote{In particular
${B}_0=+1$, ${B}_1=-\frac{1}{2}$, ${B}_2=+\frac{1}{6}$,
${B}_3=0$, ${B}_4=-\frac{1}{30}$, ${B}_5=0$, ${B}_6=+\frac{1}{42}$, ${B}_7=0$,
${B}_8=-\frac{1}{30}$, ${B}_9=0$, ${B}_{10}=+\frac{5}{66}$,
${B}_{11}=0$, ${B}_{12}=-\frac{691}{2730}$, \emph{etc}.,
see \cite[Tab.~23.2, p.~810]{abramowitz_01}, \cite[p.~5]{krylov_01}
or \cite[p.~258]{gelfond_01} for further values. Note also
that some authors may use slightly different definitions for the
Bernoulli numbers, see e.g.~\cite[p.~91]{hagen_01},
\cite[pp.~32, 71]{lindelof_01}, \cite[p.~19, \no138]{gunter_03_eng}
or \cite[pp.~3--6]{arakawa_01}.}
$H_n$ and $H^{(s)}_n$ denote the $n$th harmonic number and the $n$th
generalized harmonic number of order $s$
\begin{eqnarray}
\nonumber
H_n \,\equiv\sum_{k=1}^n \frac{1}{k}\,,\qquad\qquad\qquad
H^{(s)}_n \,\equiv\sum_{k=1}^n \frac{1}{k^s}\,,
\end{eqnarray}
respectively.
Writings $\lfloor x\rfloor$ stands for the integer part of $x$,
$\operatorname{tg}z$ for the tangent of $z$,
$\operatorname{ctg}z$ for the cotangent of $z$, $\operatorname{ch}z$
for the hyperbolic cosine of $z$, $\operatorname{sh}z$ for the
hyperbolic sine of $z$,
${\operatorname{th}}z$ for the hyperbolic tangent of $z$,
$\operatorname{cth}z$ for the hyperbolic cotangent of $z$.
In order to avoid any confusion between compositional inverse and
multiplicative inverse,
inverse trigonometric and hyperbolic functions are denoted
as $\arccos$, $\arcsin$, $\operatorname{arctg}, \ldots$ and not as
$\cos^{-1}$,
$\sin^{-1}$, $\operatorname{tg}^{-1}, \ldots{}$.
Writings $\Gamma(z)$ and $\zeta(z)$ denote respectively the gamma and
the zeta functions of argument $z$.
The Pochhammer symbol~$(z)_n$, which is also known as the generalized
factorial function, is defined as the rising factorial
$(z)_n\equiv z(z+1)(z+2)\cdots(z+n-1)=\Gamma(z+n)/\Gamma
(z)$.\footnote{For nonpositive and complex $n$, only the latter
definition $(z)_n\equiv\Gamma(z+n)/\Gamma(z)$ holds.}\up{,}\footnote{Note that some writers (mostly German-speaking)
call such a function \emph{facult\'e analytique} or \emph{Facult\"at}, see e.g.~\cite{schlomilch_04}, \cite[p.~186]{schlomilch_05},
\cite[vol.~II, p.~12]{schlomilch_06}, \cite[p.~119]{hagen_01}, \cite{kramp_01}. Other names and notations
for $(z)_n$ are briefly discussed in \cite[pp.~45--47]{jordan_01} and
in \cite[pp.~47--48]{knuth_01}.\\[-8mm]} For sufficiently large $n$, not
necessarily integer,
the latter can be given by this useful approximation\looseness=1
\be\label{lk2093mffmnjw}
\begin{array}{ll}
\displaystyle
(z)_n \;&\displaystyle =\,\frac{\,n^{n+z-\frac{1}{2}}\sqrt{2\pi} \,}{\Gamma(z)\,e^{n}} 
\left\{1+ \frac{\,6 z^2 - 6z + 1\,}{12 n} + \frac{\,36 z^4 - 120 z^3 + 120 z^2  - 36 z  + 1}{288 n^2} + O(n^{-3})\right\}\\[8mm]
&\displaystyle
\,=\,\frac{\,n^z\cdot \Gamma(n)\,}{\Gamma(z)}\left\{1+ \frac{\,z(z-1)\,}{2 n}
+ \frac{\,z(z-1)(z-2)(3z-1)\,}{24 n^2} +O(n^{-3})\right\}
\end{array}
\ee
which follows from the Stirling formula for the $\Gamma
$-function.\footnote{A simpler
variant of the above formula may be found in \cite{tricomi_01}.}
Unsigned (or signless) and signed Stirling numbers of the first kind,
which are also known as \emph{factorial coefficients},
are denoted as $|S_1(n,l)|$ and $S_1(n,l)$ respectively (the latter are
related to the former
as $S_1(n,l)=(-1)^{n\pm l}|S_1(n,l)|$).\footnote{There exist more than
50 notations for the Stirling numbers,
see e.g.~\cite{gould_02}, \cite[pp.~vii--viii, 142, 168]{jordan_01},
\cite[pp.~410--422]{knuth_02}, \cite[Sect.~6.1]{knuth_01}, and we do
not insist on our particular notation, which may
seem for certain not properly chosen.} Because in literature various
names, notations and definitions were adopted
for the Stirling numbers of the first kind, we specify that
we use exactly the same definitions and notation as in \cite[Section
2.1]{iaroslav_08},
that is to say $|S_1(n,l)|$ and $S_1(n,l)$ are defined as the coefficients
in the expansion of rising/falling factorial
\be\label{x2l3dkkk03d}
\specialnumber{a,b}
\begin{cases}
\displaystyle 
\prod_{k=0}^{n-1} (z+k) 
\,=\,(z)_n\,=\,\frac{\Gamma(z+n)}{\Gamma(z)}\,=\,\sum_{l=1}^n |S_1(n,l)|\cdot z^l \,=\, \sum_{l=0}^\infty |S_1(n,l)|\cdot z^l\\[6mm]
\displaystyle 
\prod_{k=0}^{n-1} (z-k) 
\,=\,(z-n+1)_n\,=\frac{\Gamma(z+1)}{\Gamma(z+1-n)}\,=\,\sum_{l=1}^n S_1(n,l)\cdot z^l \,=\,\sum_{l=0}^\infty S_1(n,l)\cdot z^l 
\end{cases}
\ee
respectively, where $z\in\mathbb{C}$ and $n\geqslant1$. Note that
if $l\notin[1,n]$, where $l$ is supposed to be
nonnegative, then $S_1(n,l)=0$, except for $S_1(0,0)$ which is set to 1
by convention.
Alternatively, the same numbers may be equally defined as the
coefficients in the following MacLaurin series
\be\label{ld2jr3mnfdmd}
\specialnumber{a,b}
\begin{cases}
\displaystyle
(-1)^l\frac{\ln^l(1-z)}{l!}\,=\sum_{n=l}^\infty\!\frac{|S_1(n,l)|}{n!}z^n \,=\sum_{n=0}^\infty\!\frac{|S_1(n,l)|}{n!}z^n \,, \qquad & |z|<1\,,\quad l=0, 1, 2, \ldots \\[6mm]
\displaystyle
\frac{\ln^l(1+z)}{l!}\,=\sum_{n=l}^\infty\!\frac{S_1(n,l)}{n!}z^n \,=\sum_{n=0}^\infty\!\frac{S_1(n,l)}{n!}z^n\,, \qquad & |z|<1\,,\quad l=0, 1, 2, \ldots
\end{cases}
\ee
Signed Stirling numbers of the first kind, as we defined them above,
may be also given via the following explicit formula
\be\label{io20323m3e}
S_1(n,l)\,=\,
\frac{(2n-l)!}{(l-1)!}
\sum_{k=0}^{n-l}\frac{1}{(n+k)(n-l-k)!(n-l+k)!}
\sum_{r=0}^{k}\frac{(-1)^{r} r^{n-l+k} }{r!(k-r)!} 
\ee
$l\in[1,n]$,
which may be useful for the computation of $S_1(n,l)$ when $n$ is not
very large.\footnote{From the above definitions, it follows that:
$S_1(1,1)=+1$, $S_1(2,1)=-1$, $S_1(2,2)=+1$, $S_1(3,1)=+2$,
$S_1(3,2)=-3$, $S_1(3,3)=+1$, \ldots\,,
$S_1(8,5)=-1960$, \ldots\,, $S_1(9,3)=+118\,124$, \emph{etc.} Note
that there is an error in Stirling's treatise \cite{stirling_01}:
in the last line in the table on p.~11 \cite{stirling_01} the value of
$|S_1(9,3)|=118\,124$ and not 105\,056. This error has been
noted by Jacques Binet \cite[p.~231]{binet_01}, Charles Tweedie \cite[p.~10]{tweedie_01} and some others (it was also corrected
in some translations of \cite{stirling_01}).}
All three above definitions agree with those adopted by Jordan \cite[Chapt.~IV]{jordan_01}, \cite{jordan_02,jordan_00}, Riordan
\cite[p.~70 \emph{et seq.}]{riordan_01},
Mitrinovi\'c \cite{mitrinovic_01}, Abramowitz \& Stegun \cite[\no
24.1.3, p.~824]{abramowitz_01} and many others (moreover,
modern CAS, such as \textsl{Maple} or \textsl{Mathematica}, also share
these definitions; in particular \texttt{Stirling1(n,l)} in the former
and \texttt{StirlingS1[n,l]} in the latter
correspond to our $S_1(n,l)$).\footnote{A quick analysis of several
alternative names, notations and definitions may be found in works of
Charles Jordan \cite[pp.~vii--viii, 1 and Chapt.~IV]{jordan_01}, Gould
\cite{gould_02,gould_03}, and Donald E.~Knuth
\cite[Sect.~6.1]{knuth_01}, \cite[pp.~410--422]{knuth_02}.\label{alkjcow2edchb}}
Kronecker symbol (or Kronecker delta) of arguments $l$ and $k$ is denoted
by $\,\delta_{l,k}\,$ ($\,\delta_{l,k}=1\,$
if $l=k$ and $\,\delta_{l,k}=0\,$ if $l\neq k$).
$\operatorname{Re}{z}$ and $\operatorname{Im}{z}$ denote respectively
real and imaginary parts
of $z$. 
Letter $i$ is never used as index and is $\sqrt{-1\,}$. The writing
$\operatorname{res}_{z=a} f(z)$ stands for
the residue of the function $f(z)$
at the point $z=a$. Finally, by the relative error between the quantity $A$
and its approximated value $B$, we mean $(A-B)/A$.
Other notations are standard.\looseness=-1

\section{A convergent series representation for generalized Euler's constants $\gamma_m$
involving Stirling numbers and polynomials in $\pi^{-2}$}
\subsection{Derivation of the series expansion}
In 1893 Johan Jensen \cite{jensen_04,jensen_03} by contour
integration methods
obtained an integral formula for the $\zeta$-function
\be\label{kjd02jddnsa}
\begin{array}{cc}
\displaystyle
\zeta(s) = \frac{1}{s-1} + \frac{1}{2} + 2\!\!\int\limits_0^{\pi/2} \! 
\frac{(\cos\theta)^{s-2}\sin s\theta}{e^{2\pi\tg\theta}-1} d\theta  \,
=\,
\frac{1}{s-1} + \frac{1}{2} + 2\!
\int\limits_0^\infty \! 
\frac{\sin(s \arctg x)\,}{\left(e^{2\pi x}-1\right) \left(x^2+1\right)^{s/2}}\, dx \\[8mm]
\displaystyle
=\,
\frac{1}{s-1} + \frac{1}{2} + \frac{1}{i}\!
\int\limits_0^\infty \! 
\frac{(1-ix)^{-s}-(1+ix)^{-s}}{\,e^{2\pi x}-1\,} \, dx \,,\qquad\quad s\neq 1
\end{array}
\ee
which extends {\eqref{kj023dndndr3}} to the entire complex plane except
$s=1$. Expanding the above formula
into the Laurent series about $s=1$ and performing the term-by-term
comparison of the
derived expansion with the Laurent series {\eqref{dhd73vj6s1}} yields
the following representation for the $m$th Stieltjes constant
\be\label{kljc3094jcmfd}
\gamma_m \,=\,\frac{1}{2}\delta_{m,0}+\,\frac{1}{i}\!\int\limits_0^\infty \! \frac{dx}{e^{2\pi x}-1} \left\{
\frac{\ln^m(1-ix)}{1-ix} - \frac{\ln^m(1+ix)}{1+ix} 
\right\}\,,
\qquad m=0, 1, 2,\ldots
\ee
which is due to the Jensen and Franel.\footnote{In the explicit
form, this integral formula was given by Franel in 1895 \cite{franel_01} (in the above, we corrected the original Franel's formula which was not valid for $m=0$).
However, it was remarked by Jensen
\cite{jensen_03} that it can be elementary derived from {\eqref{kjd02jddnsa}} obtained two years earlier and it is hard to disagree
with him.
By the way, it is curious that in works of modern authors, see
e.g.~\cite{connon_01,choi_01}, formula {\eqref{kljc3094jcmfd}}
is often attributed to Ainsworth and Howell, who discovered it
independently much later \cite{ainsworth_01}.}
Making a change of variable in the latter formula \break \mbox{$\,x=-\frac{1}{2\pi}\ln(1-u)\,$}, we have
\be\label{jhvc94hfhnf}
\gamma_m \,=\,\frac{1}{2}\delta_{m,0}+\,\frac{1}{\,2\pi i\,}\!\!\!\bigints\limits_{\!\!\!\!\!\!\!\!\!0}^{\;\;\;\;\;\;1} \!\! \left\{
\dfrac{\ln^m\!\left[1-\dfrac{\ln(1-u)}{2\pi i}\right]}{\,1-\dfrac{\ln(1-u)}{2\pi i}\,} \,-\, 
\dfrac{\ln^m\!\left[1+\dfrac{\ln(1-u)}{2\pi i}\right]}{\,1+\dfrac{\ln(1-u)}{2\pi i}\,} 
\right\}\frac{du}{\,u\,} 
\ee
where $\, m=0, 1, 2,\ldots$

Now, in what follows, we will use a number of basic properties of
Stirling numbers, which can be found in an amount sufficient for the
present purpose
in the following literature: \cite{stirling_01,hindenburg_01,kramp_01}, \cite[Book I, part I]{laplace_02}, \cite{ettingshausen_01,schlaffli_01,schlaffli_02,schlomilch_04},
\cite[pp.~186--187]{schlomilch_05}, \cite[vol.~II,
pp.~23--31]{schlomilch_06}, \cite{appel_01,cayley_00,cayley_01,cayley_02,boole_01,glaisher_02}, \cite[p.~129]{carlitz_02},
\cite[Chapt.~IV]{jordan_01}, \cite{jordan_02,jordan_00,nielsen_04}, \cite[pp.~67--78]{nielsen_01}, \cite{nielsen_03,tweedie_01},
\cite[Sect.~6.1]{knuth_01}, \cite[pp.~410--422]{knuth_02}, \cite[Chapt.~V]{comtet_01}, \cite{dingle_01},
\cite[Chapt.~4, \S3, \no196--\no210]{polya_01_eng}, \cite[p.~60 \emph
{et seq.}]{hagen_01}, \cite{netto_01},
\cite[p.~70 \emph{et seq.}]{riordan_01}, \cite[vol.~1]{stanley_01},
\cite{bender_01}, \cite[Chapt.~8]{charalambides_01}, \cite[\no
24.1.3, p.~824]{abramowitz_01}, \cite[Sect.~21.5-1, p.~824]{korn_01},
\cite[vol.~III, p.~257]{bateman_01}, \cite{norlund_02,steffensen_02}, \cite[pp.~91--94]{conway_01}, \cite[pp.~2862--2865]{weisstein_04},
\cite[Chapt.~2]{arakawa_01}, \cite{mitrinovic_01,gould_01,gould_02,gould_03,wachs_01,carlitz_02,carlitz_03}, \cite[p.~642]{olson_01}, \cite{salmieri_01,gessel_01,wilf_01,moser_01,bellavista_01,wilf_02,temme_02,howard_01,butzer_02,butzer_01,hwang_01,adamchik_03,timashev_01,grunberg_01,louchard_01,shen_01,shirai_01,sato_01,rubinstein_01,rubinstein_02,hauss_01,skramer_01,iaroslav_08}. Note that many writers discovered
these numbers independently, without realizing that they deal with the
Stirling numbers.
For this reason, in many sources, these numbers may appear under
different names, different notations and
even slightly different definitions.\footnote{Actually, only in the
beginning of the XXth century, the name ``Stirling numbers'' appeared
in mathematical literature
(mainly, thanks to Thorvald N.~Thiele
and Niels Nielsen \cite{nielsen_04,tweedie_01}, \cite[p.~416]{knuth_02}).
Other names for these numbers include: \emph{factorial coefficients},
\emph{faculty's coefficients} (\emph{Facult\"atencoefficienten},
\emph{coefficients de la facult\'e analytique}), \emph{differences of
zero} and even \emph{differential coefficients of nothing}.
Moreover, the Stirling numbers are also closely connected to the \emph
{generalized Bernoulli numbers} $B^{(s)}_n$, also known as
\emph{Bernoulli numbers of higher order}, see e.g.~\cite[p.~129]{carlitz_02}, \cite[p.~449]{gould_01}, \cite[p.~116]{gould_02};
many of their properties may be, therefore, deduced from those of $B^{(s)}_n$.}

\begin{figure}[!t]   
\centering
\includegraphics[width=0.8\textwidth]{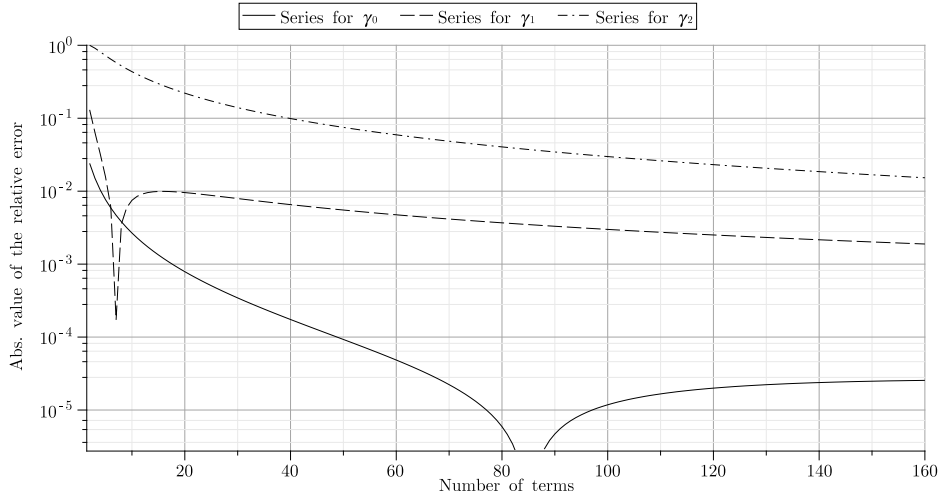}\vspace{-14mm}
\caption{Absolute values of relative errors of the series expansion for $\,\gamma_0\,$, $\,\gamma_1\,$ and 
$\,\gamma_2\,$ given by \protect\eqref{jkhf3984fhd}--\protect\eqref{jcwio0ecn32}, logarithmic scale.}
\label{kjfc0234nd}
\end{figure}

Consider the generating equation for the unsigned Stirling numbers of
the first kind, formula (\ref{ld2jr3mnfdmd}a).
This power series is uniformly and absolutely convergent inside the
disk $|z|<1$.
Putting $l+m-1$ instead of $l$, multiplying both sides by $(l)_m$ and
summing over $l=[1,\infty)$,
we obtain for the left side
\be\notag
\begin{array}{ll}
\displaystyle
\sum_{l=1}^\infty   (l)_m \cdot \frac{\big[-\ln(1-z)\big]^{l+m-1}}{(l+m-1)!}
\,=\,\sum_{l=1}^\infty \frac{\big[-\ln(1-z)\big]^{l+m-1}}{(l-1)!} =\\[6mm]
\displaystyle\qquad\qquad
=\,\big[-\ln(1-z)\big]^m \!\cdot\underbrace{\sum_{l=1}^\infty \frac{\big[-\ln(1-z)\big]^{l-1}}{(l-1)!}}_{e^{-\ln(1-z)}} \,
=\, (-1)^m\cdot\frac{\ln^m(1-z)}{1-z}
\end{array}
\ee
while the right side of (\ref{ld2jr3mnfdmd}a), in virtue of the
absolute convergence, becomes
\begin{eqnarray*}
\displaystyle
\sum_{l=1}^\infty\,(l)_m\!\cdot \sum_{n=0}^\infty \frac{\big
|S_1(n,l+m-1)\big|}{n!} \, z^n
& =&
\sum_{n=0}^\infty\frac{z^n}{n!} \,\cdot\!\!\!
\sum_{l=1}^{n-m+1}\!\!(l)_m\!\cdot\big|S_1(n,l+m-1)\big| \\[7mm]
&= & m!\cdot\!
\sum_{n=0}^\infty\frac{\,\big|S_1(n+1, m+1)\big|\,}{n!} \,z^n
\end{eqnarray*}
Whence
%
%e13 #&#
\begin{eqnarray}
\label{iu2d092n1}
\frac{\ln^m(1-z)}{1-z}\,=\,(-1)^m m!\cdot\!
\sum_{n=0}^\infty\frac{\,\big|S_1(n+1, m+1)\big|\,}{n!} \,z^n\,,
\qquad
\begin{array}{l}
m=0, 1, 2,\ldots \\[6pt]
|z|<1
\end{array}
\end{eqnarray}
Writing in the latter $-z$ for $z$,
and then subtracting
one from another yields the following series
%
%e14 #&#
\be
\label{djhd9ehdbne}
\frac{\ln^m(1-z)}{1-z}-\frac{\ln^m(1+z)}{1+z}\, =\,2(-1)^m m!\cdot\!
\sum_{k=0}^\infty\frac{\,\big|S_1(2k+2,m+1)\big|\,}{(2k+1)!}
\,z^{2k+1}
\qquad\quad
\ee
$m=0, 1, 2,\ldots\,$,
which is absolutely and uniformly convergent in the unit disk $|z|<1$, and
whose coefficients grow logarithmically with $k$
%
%e15 #&#
\begin{eqnarray}
\label{j6s8g64r}
\frac{\,\big|S_1(2k+2,m+1)\big|\,}{(2k+1)!}\sim\frac{\,\ln^m{k}\,
}{m!}\,,\qquad\quad k\to\infty\,,\qquad m=0, 1, 2,\ldots
\end{eqnarray}
in virtue of known asymptotics for the Stirling numbers,
see e.g.~\cite[p.~261]{jordan_00}, \cite[p.~161]{jordan_01}, \cite[\no24.1.3, p.~824]{abramowitz_01}, \cite[p.~348, Eq.~(8)]{wilf_02}.
Using formul{\ae}~from \cite[p.~217]{comtet_01}, \cite[p.~1395]{shen_01}, \cite[p.~425, Eq.~(43)]{kowalenko_01},
the law for the formation of first coefficients may be also written in
a more simple form
\be\label{uf87tfuy89}
\frac{\,\big|S_1(2k+2,m+1)\big|\,}{(2k+1)!}\,=\,
\begin{cases}
\,1 \,,  & m=0\\[1mm]
\,H_{2k+1}\,,& m=1 \\[1mm]
\,\frac{1}{2}\big\{H^2_{2k+1} - H^{(2)}_{2k+1}\big\} \,,\qquad\quad& m=2 \\[1mm]
\,\frac{1}{6}\big\{H^3_{2k+1} - 3H_{2k+1} H^{(2)}_{2k+1}+2H^{(3)}_{2k+1}\big\}\,,\qquad\quad& m=3   %\\[1mm]
\end{cases}
\ee
For higher $m$, values of this coefficient may be similarly reduced to
a non-linear combination of the generalized harmonic numbers.
Since expansion {\eqref{djhd9ehdbne}} holds only inside the unit circle,
it cannot be directly used for the insertion into Jensen--Franel's
integral formula {\eqref{kljc3094jcmfd}}.
However, if we put in {\eqref{djhd9ehdbne}} $z=\frac{1}{2\pi i}\ln
(1-u)$, we obtain for the right part
\be\label{jwoerivh304}
\begin{array}{ll}
\displaystyle
2\,(-1)^m m!\cdot\!\sum_{k=0}^\infty\frac{\,\big|S_1(2k+2,m+1)\big|\,}{\,(2\pi i)^{2k+1}}\cdot\underbrace{\frac{\ln^{2k+1}(1-u)}{(2k+1)!}
}_{\text{see (\ref{ld2jr3mnfdmd}a)}} = \\[10mm]
\displaystyle\qquad\qquad
\,=\,2i\,(-1)^m m!\cdot\!\sum_{k=0}^\infty\frac{\,(-1)^k \big|S_1(2k+2,m+1)\big|\,}{\,(2\pi)^{2k+1}}\cdot\!\sum_{n=1}^\infty\!\frac{\big|S_1(n,2k+1)\big|}{n!}\,u^n \\[7mm]
\displaystyle\qquad\qquad
=\,2i\,(-1)^m m!\cdot\!\sum_{n=1}^\infty\!\frac{\,u^n\,}{n!} \cdot\!\sum_{k=0}^{\lfloor\!\frac{1}{2}n\!\rfloor}
\frac{\,(-1)^k \big|S_1(2k+2,m+1)\big|\cdot\big|S_1(n,2k+1)\big|\,}{\,(2\pi)^{2k+1}}
\end{array}
\ee
Therefore, for $m=0, 1, 2,\ldots$\,, we have
\be\label{joc20j4cxcs}
\begin{array}{ll}
\displaystyle
\frac{1}{\,2\pi i\,}
\left\{
\dfrac{\ln^m\!\left[1-\dfrac{\ln(1-u)}{2\pi i}\right]}{\,1-\dfrac{\ln(1-u)}{2\pi i}\,} \,-\, 
\dfrac{\ln^m\!\left[1+\dfrac{\ln(1-u)}{2\pi i}\right]}{\,1+\dfrac{\ln(1-u)}{2\pi i}\,} 
\right\}=\, \\[10mm]
\displaystyle\qquad\qquad\qquad=\,
\frac{\,(-1)^m m!\,}{\pi}\!\sum_{n=1}^\infty\!\frac{\,u^n\,}{n!} 
\cdot\!\sum_{k=0}^{\lfloor\!\frac{1}{2}n\!\rfloor}\frac{\,(-1)^k \big|S_1(2k+2,m+1)\big| \cdot\big|S_1(n,2k+1)\big|\,}{\,(2\pi)^{2k+1}}
\end{array}
\ee
which uniformly holds in $|u|<1$ and also is valid for $u=1$.\footnote
{The unit radius of convergence of this series is conditioned
by the singularity the most closest to the origin. Such singularity is
a branch point located at $u=1$. Note also that since the series
is convergent for $u=1$ as well, in virtue of Abel's theorem on power
series, it is uniformly convergent everywhere on the disc $|u|\leqslant
1-\varepsilon$,
where positive parameter $\varepsilon$ can be made as small as we please.}
Substituting {\eqref{joc20j4cxcs}} into {\eqref{jhvc94hfhnf}} and
performing the term-by-term integration from $u=0$ to $u=1$ yields
the following series representation for $m$th generalized Euler's constant
\be\label{jkhf3984fhd}
\gamma_m\,=\,\frac{1}{2}\delta_{m,0}+
\frac{\,(-1)^m m!\,}{\pi} \!\sum_{n=1}^\infty\frac{1}{\,n\cdot n!\,} \!
% \sum_{k=0}^{\lfloor\!\frac{1}{2}n\!\rfloor}\frac{\,(-1)^{k}\!\cdot a_m(k)\cdot\big|S_1(n,2k+1)\big|}{\,(2\pi)^{2k+1}\,}\,,\qquad\quad m=0, 1, 2,\ldots
\sum_{k=0}^{\lfloor\!\frac{1}{2}n\!\rfloor}\frac{\,(-1)^{k}\big|S_1(2k+2,m+1)\big| \cdot\big|S_1(n,2k+1)\big|\,}{\,(2\pi)^{2k+1}\,}
\ee
where $m=0, 1, 2,\ldots{}$.
In particular, for Euler's constant and first Stieltjes constant, we
have following series expansions
\be\label{jcwio0ecn32}
\begin{array}{rcl}
\displaystyle
\gamma\; &=&\displaystyle\,\frac{1}{2}+
\frac{\,1\,}{\pi}\sum_{n=1}^\infty\frac{1}{\,n\cdot n!\,} 
\sum_{k=0}^{\lfloor\!\frac{1}{2}n\!\rfloor}\frac{\,(-1)^{k}\!\cdot (2k+1)!\cdot\big|S_1(n,2k+1)\big|}{\,(2\pi)^{2k+1}\,} \\[5mm]
\displaystyle
& =&\,\displaystyle \frac{1}{2}+ \frac{1}{2\pi^2}+\frac{1}{8\pi^2}+\frac{1}{18}\!\left(\frac{1}{\pi^2}-\frac{3}{4\pi^4}\right) 
+\frac{3}{96}\!\left(\frac{1}{\pi^2}-\frac{3}{2\pi^4}\right)  \\[6mm]
& &\displaystyle
+\frac{1}{600}\!\left(\frac{12}{\pi^2}-\frac{105}{4\pi^4}+\frac{15}{4\pi^6}\right)  
+\frac{1}{4\,320}\!\left(\frac{60}{\pi^2}-\frac{675}{4\pi^4}+\frac{225}{4\pi^6}\right)  + \ldots \\[7mm]
\displaystyle
\gamma_1\:& =&\displaystyle\,-\frac{\,1\,}{\pi}\sum_{n=1}^\infty\frac{1}{\,n\cdot n!\,} 
\sum_{k=0}^{\lfloor\!\frac{1}{2}n\!\rfloor}\frac{\,(-1)^{k}\!\cdot (2k+1)!\cdot H_{2k+1}\cdot\big|S_1(n,2k+1)\big|}{\,(2\pi)^{2k+1}\,} \\[5mm]
\displaystyle
& =&\,\displaystyle -\frac{1}{2\pi^2}-\frac{1}{8\pi^2}-\frac{1}{18}\!\left(\frac{1}{\pi^2}-\frac{11}{8\pi^4}\right) 
-\frac{3}{96}\!\left(\frac{1}{\pi^2}-\frac{11}{4\pi^4}\right)  \\[6mm]
& &\displaystyle
-\frac{1}{600}\!\left(\frac{12}{\pi^2}-\frac{385}{8\pi^4}+\frac{137}{16\pi^6}\right) 
-\frac{1}{4\,320}\!\left(\frac{60}{\pi^2}-\frac{2\,475}{8\pi^4}+\frac{2055}{16\pi^6}\right)  - \ldots
\end{array}
\ee
respectively.
As one can easily notice, each coefficient of these expansions contains
polynomials in $\pi^{-2}$ with rational coefficients.
The rate of convergence of this series, depicted in {Fig.~\ref{kjfc0234nd}}, is relatively slow and depends, at least for the moderate
number of terms, on~$m$:
the greater the order $m$, the slower the convergence.
A more accurate description of this dependence, as well as the exact
value of the rate of convergence,
both require a detailed convergence analysis of {\eqref{jkhf3984fhd}},
which is performed in the next section.
%%%%%%%%%%%%%%%%%%%%%%%%%%%%%%%%%%%%%%%%%%%%%%%%%%%%%%%%%%
%%%%%%%%%%%%%%%%%%%%%%%%%%%%%%%%%%%%%%%%%%%%%%%%%%%%%%%%%

%s2.2 #&#
\subsection{Convergence analysis of the derived series}
The convergence analysis of series {\eqref{jkhf3984fhd}} consists in the
study of its general term, which is given by the finite truncated sum
over index $k$.
This sum has only odd terms, and hence, by elementary transformations,
may be reduced
to that containing both odd and even terms
\be\label{lkh908gb9gi8vityr}
\begin{array}{ll}
& \displaystyle 
\sum_{k=0}^{\lfloor\!\frac{1}{2}n\!\rfloor}(-1)^{k}\frac{\,\big|S_1(2k+2,m+1)\big| \cdot\big|S_1(n,2k+1)\big|\,}{\,(2\pi)^{2k+1}\,} \,=\\[6mm]
& \displaystyle \qquad\qquad
= \sum_{k=0}^{\lfloor\!\frac{1}{2}n\!\rfloor} (-1)^{\frac{1}{2}(2k+1)-\frac{1}{2}}
\frac{\,\big|S_1(2k+1+1,m+1)\big| \cdot\big|S_1(n,2k+1)\big|\,}{\,(2\pi)^{2k+1}\,}\\[8mm]
& \displaystyle \qquad\qquad
= \,\frac{1}{2}\!\sum_{l=1}^{n} \big[1-(-1)^{l}\big] \cdot(-1)^{\frac{1}{2}(l-1)}\cdot
\frac{\,\big|S_1(l+1,m+1)\big| \cdot\big|S_1(n,l)\big|\,}{\,(2\pi)^{l}\,}
\,=\,\ldots
\end{array}
\ee
where, in the last sum, we changed the summation index by putting $l=2k+1$.
Now, from the second integral formula for the unsigned Stirling numbers
of the first kind,
see {\eqref{ock2w3jkmd1}}, it follows that
\be\notag
\begin{array}{ll}
\displaystyle
(-1)^{\frac{1}{2}(l-1)}\cdot\frac{\,\big|S_1(l+1,m+1)\big|\,}{\,(2\pi)^{l}\,}\,
=\,\frac{(-1)^m}{\,2\pi \,}\cdot\frac{(l+1)!}{(m+1)!}\!\!
\ointctrclockwise\limits_{|z|=r}\!\!\left[+\frac{i}{2\pi z}\right]^l\frac{\ln^{m+1}(1-z)}{z^2}\, dz\,\\[10mm]
\displaystyle
(-1)^l\cdot(-1)^{\frac{1}{2}(l-1)}\cdot\frac{\,\big|S_1(l+1,m+1)\big|\,}{\,(2\pi)^{l}\,}\,
=\,\frac{(-1)^m}{\,2\pi \,}\cdot\frac{(l+1)!}{(m+1)!}\!\!
\ointctrclockwise\limits_{|z|=r}\!\!\left[-\frac{i}{2\pi z}\right]^l\frac{\ln^{m+1}(1-z)}{z^2}\, dz
\end{array}
\ee
where $0<r<1$. Therefore, since $(l+1)!=\int \! x^{l+1} e^{-x} dx$
taken from $0$ to $\infty$,
the last sum in {\eqref{lkh908gb9gi8vityr}} reduces to the following
integral representation
\be\label{oi23jrn3ds3}
\begin{array}{ll}
\ldots\;&\displaystyle  =\,
\frac{(-1)^m}{\,4\pi (m+1)! \,}\sum_{l=1}^{n} \big|S_1(n,l)\big|\cdot(l+1)!\cdot\!\! \ointctrclockwise\limits_{|z|=r}\!\!
\left[\left(\frac{i}{2\pi z}\right)^{\!l} - \left(-\frac{i}{2\pi z}\right)^{\!l}\right]\frac{\ln^{m+1}(1-z)}{z^2}\, dz \\[9mm]
&\displaystyle 
=\,\frac{(-1)^m}{\,4\pi (m+1)! \,}\cdot\int\limits_{0}^\infty\left[\sum_{l=1}^{n} \big|S_1(n,l)\big|\!\! \ointctrclockwise\limits_{|z|=r}\!\!
\left[\left(\frac{ix}{2\pi z}\right)^{\!l} - \left(-\frac{ix}{2\pi z}\right)^{\!l}\right]\frac{\ln^{m+1}(1-z)}{z^2}\, dz\right] x\,e^{-x}\,\,dx \\[10mm]
&\displaystyle 
=\,\frac{(-1)^m}{\,4\pi (m+1)! \,}\cdot \!\! \ointctrclockwise\limits_{|z|=r}\!\!
\frac{\ln^{m+1}(1-z)}{z^2}\left\{\int\limits_{0}^\infty\left[\left(\frac{ix}{2\pi z}\right)_{\!\!n}
- \left(-\frac{ix}{2\pi z}\right)_{\!\!n}\right]x\,e^{-x}\,dx\right\} dz %\\[8mm]
\end{array}
\ee
The integral in curly brackets is difficult to evaluate in a closed-form,
but at large $n$, its asymptotical value may be readily obtained.

Function $1/\Gamma(z)$ is analytic on the entire complex $z$-plane,
and hence, can be expanded
into the MacLaurin series
\be\label{poi2d293dm}
\frac{1}{\Gamma(z)}\,=\,z+\gamma z^2+ \left(\! \frac{\gamma^2}{2}-\frac{\pi^2}{12}\right)\!z^3 +
\ldots\,\equiv\sum_{k=1}^\infty z^k a_k \,,\qquad |z|<\infty\,,\\[6mm]
\ee
where
\be\notag
a_k\,\equiv\,\frac{1}{k!}\cdot \left[\frac{1}{\Gamma(z)}\right]^{(k)}_{z=0} \!=\, 
\frac{(-1)^k}{\,\pi \, k!\,}\cdot \Big[\sin\pi x\cdot\Gamma(x)\Big]^{(k)}_{x=1} 
\ee 
see e.g.~\cite[p.~256, \no6.1.34]{abramowitz_01}, \cite[pp.~344 \&
349]{wilf_02}, \cite{hayman_01}.
Using Stirling's approximation for the Pochhammer symbol {\eqref{lk2093mffmnjw}}, we have for
sufficiently large $n$
\be\label{poi2d293dm2}
\begin{array}{l}
\displaystyle
\left(\frac{ix}{2\pi z}\right)_{\!\!n} - \left(-\frac{ix}{2\pi z}\right)_{\!\!n}\sim
\, \frac{n^{\frac{ix}{2\pi z}}\cdot\Gamma(n)}{\Gamma\left(\frac{ix}{2\pi z}\right)}\,-\,
\frac{n^{-\frac{ix}{2\pi z}}\cdot\Gamma(n)}{\Gamma\left(-\frac{ix}{2\pi z}\right)}\,=\\[8mm]
\displaystyle
=\,(n-1)! \left[\exp\!\left(\frac{ix\ln n}{2\pi z}\right)\!\sum_{k=1}^\infty a_k \!\left(\frac{i\,x}{2\pi z}\right)^{\!k}- 
\exp\!\left(-\frac{ix\ln n}{2\pi z}\right)\!\sum_{k=1}^\infty (-1)^k a_k \!\left(\frac{i\,x}{2\pi z}\right)^{\!k}\right]
\end{array}
\ee 
Substituting this approximation into the integral in curly brackets
from {\eqref{oi23jrn3ds3}},
performing the term-by-term integration\footnote{Series
(23)--(24) being
uniformly convergent.} and taking into account that $z^{-s}\Gamma
(s)=\int\! x^{s-1} e^{-zx} dx$ taken over $x\in[0,\infty)$,
yields
\be\label{c293mned}
\begin{array}{ll}
\displaystyle
\int\limits_{0}^\infty\left[\left(\frac{ix}{2\pi z}\right)_{\!n} - \left(-\frac{ix}{2\pi z}\right)_{\!n}\right]x\,e^{-x}\,dx\,\sim  \\[7.8mm]
\displaystyle\qquad\qquad
\sim \,(n-1)!\sum_{k=1}^\infty a_k \!\left(\frac{i}{2\pi z}\right)^{\!k}\!\!\cdot (k+1)! \left\{\left[1+\frac{i\ln n}{2\pi z}\right]^{-k-2} \!\!\! -
(-1)^k\left[1-\frac{i\ln n}{2\pi z}\right]^{-k-2} \right\}\\[8mm]
\displaystyle\qquad\qquad
\sim \,(n-1)! \cdot\frac{32\,i\pi^3 z^3 \big(4\pi^2 z^2 - 3\ln^2 n\big)}{\,\big(4\pi^2 z^2 +\ln^2 n\big)^3\,}\,,
\qquad\qquad\quad n\to\infty\,,
\end{array}
\ee
where, at the final stage, we retained only the first significant term
corresponding to factor $k=1$.\footnote{The
second term of this sum, corresponding to $k=2$, is
\begin{eqnarray*}
(n-1)!\cdot \frac{396\,i\gamma\pi^3 z^3 
\big(4\pi^2 z^2 - \ln^2 n\big)\ln n}{\big(4\pi^2 z^2 +\ln^2 n\big)^{4}}=(n-1)!
\cdot O\left(\frac{1}{\,\ln^{5} n\,}\right)\,,\qquad\quad n\to\infty\,,
\end{eqnarray*}
and hence, may be neglected at large $n$.}
Now, if $|z|\leqslant1-e^{-1}\approx0.63$, then the principal
branch of \mbox{$\big|\ln^{m+1}(1-z)\big|\leqslant 1$}
independently of $m$ and $\arg{z}$. Analogously, one can always find
such sufficiently large~$n_0$, that for any however small $\varepsilon>0$,
\be\label{29ci23nd3jw}
\left| \frac{32\,i\pi^3 z^3 \big(4\pi^2 z^2 - 3\ln^2 n\big) }{\,\big(4\pi^2 z^2 +\ln^2 n\big)^3\,}\right|<\varepsilon\,,
\qquad n\geqslant n_0\,,
\ee
on the circle $|z|=1-e^{-1}$
(for example, if $\varepsilon=1$, then $n_0=1222$;
if $\varepsilon=0.1$, then $n_0=38\,597$; if $\varepsilon=0.01$, then
$n_0=33\,220\,487$; \emph{etc.}).\footnote{Note that
for fixed $n$, the left-hand side of {\eqref{29ci23nd3jw}} reaches its
maximum when $z$ is imaginary pure.}
Combining all these results and taking into account that $|dz|=|z|\,
d\arg{z}$, we conclude that
\be\label{k039dm3dmedc}
\begin{array}{ll}
\displaystyle
\frac{1}{\,n\cdot n!\,} \left|
\sum_{k=0}^{\lfloor\!\frac{1}{2}n\!\rfloor}\frac{\,(-1)^{k}
\big|S_1(2k+2,m+1)\big| \cdot\big|S_1(n,2k+1)\big|\,}{\,(2\pi)^{2k+1}\,} \right| <
% \, < \, \frac{1}{\,n^2\,}\cdot\frac{\varepsilon}{\,2\big(1-e^{-1}\big) (m+1)!\,}\,,
\\[8mm]
\displaystyle\qquad\qquad\qquad\qquad\qquad
< \, \frac{1}{\,n^2\,}\cdot\frac{\varepsilon}{\,2\big(1-e^{-1}\big) (m+1)!\,}
\,<\,\frac{C}{\,n^2\,}\,,
\qquad
\qquad n\geqslant n_0\,,
\end{array}
\ee
Numerical simulations, 
%%%%%%%%%%%%%%%%%%%%%%%%%%%%%%%%%%%%%%%%%%%%%%%%%%%%%%%%%%%%
\begin{figure}[!t]   
\centering
\includegraphics[width=0.8\textwidth]{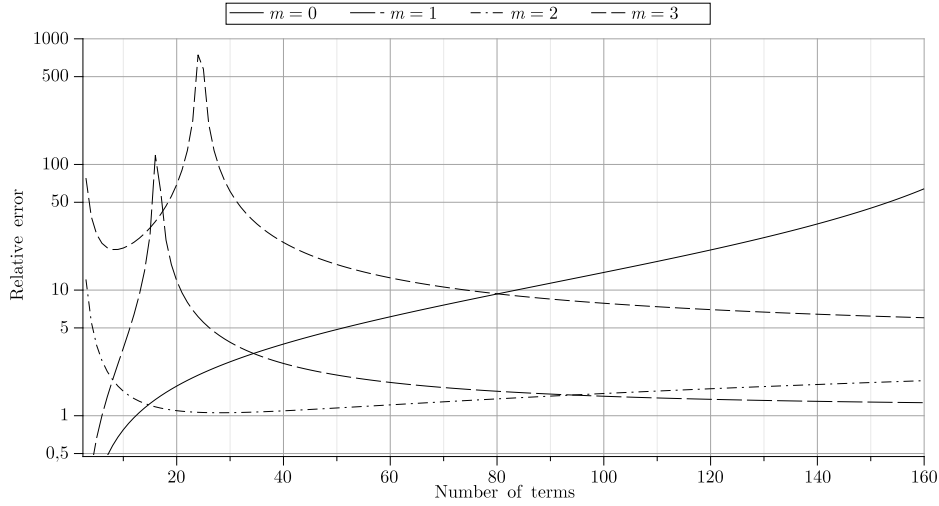}\vspace{-3mm}
\caption{Relative error between the upper bound and the left--hand side in \eqref{k039dm3dmedc} as a function of $n$ 
for four different orders $m$, logarithmic scale (the curve in long dashes correspond to $m=1$, that in short ones to $m=3$). Results displayed above correspond to $C=1/(2\pi)$.}
\label{cj2394chcned}
\end{figure}
%%%%%%%%%%%%%%%%%%%%%%%%%%%%%%%%%%%%%%%%%%%%%%%%%%%%%%%%%%%
see Fig.~\ref{cj2394chcned}, show that this
simple inequality, valid for all $m$,
may provide more or less accurate approximation\vadjust{\eject} for the general term of {\eqref{oi23jrn3ds3}}, and this greatly depends on $m$.
Moreover, the joint analysis of {Figs.~\ref{kjfc0234nd} and \ref{cj2394chcned}} also indicates that first partial sums of
series {\eqref{jkhf3984fhd}} may behave quite irregularly. One of the
reasons of such a behaviour is that
for $1\leqslant n\leqslant53$, absolute value {\eqref{29ci23nd3jw}}
increases, and it starts to decrease only
after $n=54$.\footnote{On the circle $|z|=1-e^{-1}$, absolute value {\eqref{29ci23nd3jw}} has
one of its third-order poles at $n=e^{2\pi(1-e^{-1})}\approx53.08$.
Other poles are located either below $n=1$,
e.g.~$n=e^{-2\pi(1-e^{-1})}\approx0.02$, or are complex. More
precisely, all poles of
this expression occur at $n=\big[\cos\big(2\pi(1-e^{-1})\cos
\varphi\big)
\mp i\sin\big(2\pi(1-e^{-1})\cos\varphi\big)\big]e^{\pm2\pi
(1-e^{-1})\sin\varphi}$, where $\varphi\equiv\arg z$.}
Notwithstanding, inequality {\eqref{k039dm3dmedc}} guarantees that in
all cases, the discovered series for $\gamma_m$ given by {\eqref{jkhf3984fhd}}
converges for large $n$ not worse than Euler's series $\sum n^{-2}$.

Asymptotics {\eqref{c293mned}},
as well as the rates of convergence previously obtained in \cite{iaroslav_08},
both suggest that
the exact rate of convergence of series {\eqref{jkhf3984fhd}} may also
involve logarithms.
For instance, from \cite[Sect.~3]{iaroslav_08}, it
straightforwardly follows that the rate of
convergence of this series at $m=0$ is equal to $\sum(n\ln n)^{-2}$, see also footnote \ref{gtf1a}.
Indeed, if we replace the integral in curly brackets from {\eqref{oi23jrn3ds3}} by its first-order
approximation {\eqref{c293mned}},
and then, evaluate the sum of corresponding residues at $
z_{1,2}\equiv\pm\frac{i\ln n}{2\pi}$
\vspace{2mm}
\be\notag
\begin{array}{ll}
\displaystyle
\sum_{l=1}^2\res_{z=z_l}\!\frac{z\big(4\pi^2 z^2 - 3\ln^2 n\big)\ln(1-z)}{\,\big(4\pi^2 z^2 +\ln^2 n\big)^3\,}
\,=\,\frac{\ln^2n-4\pi^2}{8\pi^2 \big(4\pi^2 +\ln^2 n\big)^2} 
 \,\sim \frac{1}{\,8\pi^2\ln^2 n\,}\\[9mm]
\displaystyle
\sum_{l=1}^2\res_{z=z_l}\!\frac{z\big(4\pi^2 z^2 - 3\ln^2 n\big)\ln^2(1-z)}{\,\big(4\pi^2 z^2 +\ln^2 n\big)^3\,}
\,=\,\frac{\ln^2n \cdot\ln\big(4\pi^2 +\ln^2 n\big)+\ldots}{8\pi^2 \big(4\pi^2 +\ln^2 n\big)^2} 
 \,\sim \frac{2\cdot\ln\ln n}{\,8\pi^2\ln^2 n\,}\\[9mm]
\displaystyle
\sum_{l=1}^2\res_{z=z_l}\!\frac{z\big(4\pi^2 z^2 - 3\ln^2 n\big)\ln^3(1-z)}{\,\big(4\pi^2 z^2 +\ln^2 n\big)^3\,}
\,=\,\frac{3\ln^2n \cdot\big[\ln^2(2\pi +i\ln n) + \ln^2(2\pi -i\ln n) \big]+\ldots}{16\pi^2 \big(4\pi^2 +\ln^2 n\big)^2} \\[5mm]
\displaystyle \qquad  \qquad \qquad \qquad \qquad \qquad \qquad \qquad
 \,\sim \frac{3\cdot\ln^2\!\ln n}{\,8\pi^2\ln^2 n\,}
\end{array}
\ee
and so on, we find that
\be\label{owjh298enc}
\begin{array}{ll}
& \displaystyle 
\frac{1}{\,n\cdot n!\,} \sum_{k=0}^{\lfloor\!\frac{1}{2}n\!\rfloor}(-1)^{k}\frac{\,\big|S_1(2k+2,m+1)\big| \cdot\big|S_1(n,2k+1)\big|\,}{\,(2\pi)^{2k+1}\,} \,\sim\\[8mm]
& \displaystyle \qquad\quad
\sim\,\frac{8\,i\,\pi^2 \, (-1)^m}{\,n^2(m+1)! \,}\! \ointctrclockwise\limits_{|z|=r}\!\!
\frac{z\big(4\pi^2 z^2 - 3\ln^2 n\big)\ln^{m+1}(1-z)}{\,\big(4\pi^2 z^2 +\ln^2 n\big)^3\,} \, dz \,
\sim\,(-1)^{m+1}\frac{\,2\pi\,}{\, \,m!\,}\cdot\frac{\ln^m\ln n}{\,n^2\ln^2 n\,}\qquad
\end{array}
\ee
in virtue of the Cauchy residue theorem. Of course, this formula is
only a rough approximation,
because poles $z_{1,2}$ belongs to the disc $|z|=r<1$ only if
$1\leqslant n\leqslant535$, while formula {\eqref{c293mned}} is
a double first-order approximation and is
accurate only for large $n$.
Furthermore, residues were also evaluated only in the first
approximation. However, obtained expression
gives an idea of what the true rate
of convergence of series {\eqref{jkhf3984fhd}} could be, and it also
explains quite well why series for higher
generalized Euler's constants converge more slowly than those for lower
generalized Euler's constants.
Moreover, this approximation agrees with the fact that {\eqref{jkhf3984fhd}} converges not worse than Euler's series, and also is
consistent with the previously derived rate of convergence for $\gamma
$ from \cite{iaroslav_08}, which was obtained by another method.

%s3 #&#
\section{Expansion of generalized Euler's constants $\gamma_m$ into
the formal series with rational coefficients}

%s3.1 #&#
\subsection{Introduction}

Expansions into the series with rational coefficients is an interesting
and challenging subject.
There exist many such representations for Euler's constant $\gamma$
and first of them date back to
the XVIIIth century. For instance, from the Stirling series for the
digamma function
at $z=1$, it straightforwardly follows that
\be\label{lkce02m}
\gamma\,=\,\frac{1}{2} + \sum_{k=1}^{N}\frac{\,{B}_{2k}}{\,2k\,} 
+\, \theta\cdot\frac{\,{B}_{2N+2}\,}{\,2(N+1)\,} \,=\,
\frac{1}{2}+ \frac{1}{12}-\frac{1}{120}+\frac{1}{252}-\frac{1}{240}+\frac{1}{132}-\frac{691}{32\,760}+ \ldots 
\ee
where $0<\theta<1$ and $N<\infty$.\footnote{This result should be
attributed to both Stirling and De Moivre, who originated
Stirling series, see \cite[p.~135]{stirling_01} and \cite{demoivre_01} respectively (for more information on Stirling series,
see also \cite[part II, Chapter VI, p.~466]{euler_02}, \cite[p.~33]{gauss_02},
\cite[p.~329]{bromwich_01}, \cite[p.~111]{norlund_02}, \cite[\S12-33]{watson_01},
\cite[\S15-05]{jeffreys_02}, \cite[p.~530]{knopp_01}, \cite[p.~1]{copson_01},
\cite{kratzer_01}, \cite[\S4.1, pp.~293--294]{olver_01}, \cite[pp.~286--288]{gelfond_01},
\cite[\no6.1.40--\no6.1.41]{abramowitz_01}, \cite{murray_01}). Curiously, Srivastava and Choi \cite[p.~6]{srivastava_03},
did not notice the trivial connection between this series and the
Stirling series for the digamma function
and erroneously credited this result to Konrad Knopp, in whose book
\cite{knopp_01} it appears, with a slightly different remainder, on p.~527
(Knopp himself never claimed the authorship of this formula).} A more
general representation of the same kind may be obtained by
Euler--MacLaurin summation
%
%e30 #&#
\begin{eqnarray}
\gamma\,=\, \sum_{k=1}^n \frac{\,1\,}{k}-\ln{n}
-\frac{\,1\,}{2n} + \sum_{k=1}^{N} \frac{\,{B}_{2k}\,}{2k\cdot n^{2k}}
+ \theta\cdot\frac{\,{B}_{2N+2}\,}{2(N+1)\cdot n^{2N+2}}\,,
\end{eqnarray}
where $0<\theta<1$, $N<\infty$ and $n$ is positive integer, see
e.g.~\cite[\no377]{gunter_03_eng}.
Two above series are \emph{semi-convergent} (or \emph{divergent
enveloping}), i.e.~they diverge as $N\to\infty$.
The first known convergent series representation for Euler's constant
with only rational terms, as far as we know, dates back to 1790 and
is due to Gregorio Fontana and Lorenzo Mascheroni
%
%e31 #&#
\begin{eqnarray}
\label{njw3uiqch}
\gamma\,=\sum_{n=1}^{\infty} \frac{\,\big|G_n\big|\,}{n}= \frac
{1}{2}+\frac{1}{24}+\frac{1}{72}+\frac{19}{2880}+\frac{3}{800}
+\frac{863}{362\,880}+ \ldots
\end{eqnarray}
where rational coefficients $G_n$, known as \emph{Gregory's
coefficients},\footnote{These coefficients are also called
\emph{(reciprocal) logarithmic numbers}, \emph{Bernoulli numbers of
the second kind},
normalized \emph{generalized Bernoulli numbers} $B_n^{(n-1)}$, \emph
{Cauchy numbers} and normalized \emph{Cauchy numbers
of the first kind} $C_{1,n}$. They were introduced by
James Gregory in 1670 in the context of area's interpolation formula
(which is known nowadays as \emph{Gregory's interpolation formula})
and were subsequently rediscovered in various contexts by many famous
mathematicians, including Gregorio Fontana, Lorenzo Mascheroni,
Pierre--Simon Laplace, Augustin--Louis Cauchy, Jacques Binet,
Ernst Schr\"oder, Oskar Schl\"omilch, Charles Hermite, Jan C.~Kluyver
and Joseph Ser
\cite[vol.~II, pp.~208--209]{rigaud_01},
\cite[vol.~1, p.~46, letter written on November 23, 1670 to John
Collins]{newton_01}, \cite[pp.~266--267, 284]{jeffreys_02},
\cite[pp.~75--78]{goldstine_01},
\cite[pp.~395--396]{chabert_01},
\cite[pp.~21--23]{mascheroni_01}, \cite[T.~IV,
pp.~205--207]{laplace_01}, \cite[pp.~53--55]{boole_01}, \cite{van_veen_01},
\cite[pp.~192--194]{goldstine_01},
\cite{lienard_01,wachs_01,schroder_01,schlomilch_03}, \cite[pp.~65, 69]{hermite_01},
\cite{kluyver_02,ser_01}.
For more information about these important coefficients, see
\cite[pp.~240--251]{norlund_02}, \cite{norlund_01}, \cite[p.~132,
Eq.~(6), p.~138]{jordan_02}, \cite[p.~258, Eq.~(14)]{jordan_00}, \cite[pp.~266--267, 277--280]{jordan_01},
\cite{nielsen_01,nielsen_03,steffensen_01}, \cite[pp.~106--107]{steffensen_02}, \cite{davis_02},
\cite[p.~190]{weisstein_04},
\cite[p.~45, \no370]{gunter_03_eng},
\cite[vol.~III, pp.~257--259]{bateman_01}, \cite{stamper_01}, \cite[p.~229]{krylov_01}, \cite[\no600, p.~87]{proskuriyakov_01_eng},
\cite[p.~216, \no75-a]{knopp_01}
\cite[pp.~293--294, \no13]{comtet_01}, \cite{carlitz_01,howard_02,young_01,adelberg_01,zhao_01,candelpergher_01},
\cite[Eq.~(3)]{mezo_01}, \cite{merlini_01,nemes_01},
\cite[pp.~128--129]{alabdulmohsin_01},
\cite[Chapt.~4]{arakawa_01}, \cite{skramer_01,iaroslav_08}.\label{jpqwcnqasd}} are given either via their
generating function
%
%e32 #&#
\begin{eqnarray}
\label{eq32}
\frac{z}{\ln(1+z)} = 1+\sum
_{n=1}^\infty G_n z^n ,\qquad|z|<1\,,
\end{eqnarray}
or explicitly
%
%e33 #&#
\begin{eqnarray}
\label{ldhd9ehn}
G_n=\frac{1}{n!} \sum\limits_{l=1}^{n} \frac{S_1(n,l)}{l+1}
=\frac{1}{n!}\!\int\limits_0^1 (x-n+1)_n\, dx=-\frac{B_n^{(n-1)}}{\,
(n-1)\,n!\,}\,=\,\frac{C_{1,n}}{n!}
\end{eqnarray}
This series was first studied by Fontana, who, however,
failed to find a constant to which it converges. Mascheroni identified
this \emph{Fontana's constant}
and showed that it equals Euler's constant \cite[pp.~21--23]{mascheroni_01}.
This series was subsequently rediscovered many times, in particular, by
Ernst Schr\"oder in 1879 \cite[p.~115, Eq.~(25a)]{schroder_01},
by Niels E.~N{\o}rlund in 1923 \cite[p.~244]{norlund_02},
by Jan C.~Kluyver in 1924 \cite{kluyver_02}, by Charles Jordan in 1929
\cite[p.~148]{jordan_02},
by Kenter in 1999 \cite{kenter_01}, by Victor Kowalenko in 2008 \cite{kowalenko_01,kowalenko_02}.
An expansion of a similar nature
\be\label{c32c324fv}
\gamma\,=\,1 - \sum_{n=1}^\infty \frac{C_{2,n}}{\,n\cdot(n+1)!\,}
=\,
1-\frac{1}{4}-\frac{5}{72}-\frac{1}{32}-\frac{251}{14\,400}-\frac
{19}{1728} -
\frac{19\,087}{2\,540\,160} - \ldots
\ee
where rational numbers $C_{2,n}$, known as \emph{Cauchy numbers of the
second kind}\footnote{These
numbers, called by some authors signless \emph{generalized Bernoulli numbers}
$|B_n^{(n)}|$ and signless \emph{N{\o}rlund numbers}, are much less
famous than
Gregory's coefficients $G_n$, but their study is also very interesting, see
\cite[pp.~150--151]{norlund_02}, \cite[p.~12]{davis_02}, \cite{norlund_01}, \cite[vol.~III,
pp.~257--259]{bateman_01},
\cite[pp.~293--294, \no13]{comtet_01}, \cite{howard_03,adelberg_01,zhao_01,qi_01,iaroslav_08}.}
\be
\begin{cases}
\displaystyle
\frac{z}{(1+z)\ln(1+z)}\,=\,1+\sum_{n=1}^\infty\!\frac{(-1)^n z^n C_{2,n}}{n!} \\[8mm]
\displaystyle
C_{2,n}\,=\sum\limits_{l=1}^{n} \frac{|S_1(n,l)|}{l+1} = \int\limits_0^1 (x)_n\, dx\,=\,|B_n^{(n)}|
\end{cases}
\ee
follows from a little-known series
for the digamma function given by Jacques Binet in 1839 \cite[p.~257,
Eq.~(81)]{binet_01}
and rediscovered later by Niels E.~N{\o}rlund in his monograph \cite[p.~244]{norlund_02}.\footnote{Strictly speaking, Binet found only
first four coefficients of the corresponding series for the digamma
function and incorrectly calculated the last coefficient (for $K(5)$ he
took $\frac{245}{3}$
instead of $\frac{245}{6}$ \cite[p.~237]{binet_01}), but otherwise
his method and derivations are correct. It is also notable that
Binet related coefficients $K(n)$ to the Stirling numbers and
provided two different ways for their computation, see \cite[Final remark]{iaroslav_08}.} 
Series \break

\begin{eqnarray}
&&  \gamma=1-\sum_{n=1}^\infty\sum
_{k=2^{n-1}}^{2^{n}-1} \frac{n}{(2k+1)(2k+2)}  \,=\,1-\sum_{n=1}^\infty\sum
_{k=2^{n}+1}^{2^{n+1}} \frac{(-1)^{k+1}n}{\,k\,}  \,=\, 1-\frac{1}{12}-\frac{43}{420} \notag\\[3mm]
&&  \phantom{\gamma=\,}
-\frac{20\,431}{240\,240}-\frac{2\,150\,797\,323\,119}{36\,100\,888\,223\,400}
- \frac{9\,020\,112\,358\,835\,722\,225\,404\,403}{236\,453\,376\,820\,564\,453\,502\,272\,320} - \ldots\label{lk2eojmwjksd} \\[-3mm]
\notag
\end{eqnarray}
was given in the first form by Niels Nielsen in 1897 \cite[Eq.~(6)]{nielsen_02},
and in the second form by Ernst Jacobsthal in 1906 \cite[Eqs.~(8)]{jacobsthal_01}.
The same series (in various forms) was independently obtained by
Addison in 1967 \cite{addison_01} and
by Gerst in 1969 \cite{gerst_01}. The famous series
%
%e37 #&#
\begin{eqnarray}
\label{lk2jd029jde}
\gamma= \sum_{n=2}^\infty
\frac{(-1)^n}{n} \lfloor\log_2{n}\rfloor\,=\,\frac{1}{2}-\frac
{1}{3}+\frac{1}{2}-
\frac{2}{5}+\frac{1}{3}-\frac{2}{7}+\ldots
\end{eqnarray}
was first given by Ernst Jacobsthal in 1906 \cite[Eqs.~(9)]{jacobsthal_01}
and subsequently rediscovered by many writers, including Giovanni Vacca
\cite{vacca_01},
H.F.~Sandham \cite{sandham_01}, D.F.~Barrow, M.S.~Klamkin, N.~Miller
\cite{barrow_01} and
Gerst \cite{gerst_01}.\footnote{It should be remarked that this
series is often incorrectly attributed to Vacca,
who only rediscovered it. This error initially
is due to Glaisher, Hardy and Kluyver, see e.g.~their works \cite{glaisher_01,hardy_03,kluyver_02}.
It was only much later that Stefan Kr\"amer \cite{skramer_01}
correctly attributed this series to Jacobsthal \cite{jacobsthal_01}.}
Series
\be\label{f413f14cx45} 
\gamma= \sum_{n=m}^\infty
\frac{ \beta_n }{n} \lfloor\log _m{n}\rfloor ,\qquad
\beta_n = 
\begin{cases}
m-1 , & n = \mbox{multiple of } m \vspace{3pt}\cr
-1 , & n \neq\mbox{multiple of } m
\end{cases}
\ee
which generalizes foregoing Jacobsthal--Vacca's series {\eqref{lk2jd029jde}},
is due to Jan C.~Kluyver who discovered it in 1924 \cite{kluyver_02}.

In contrast, as concerns generalized Euler's constants $\gamma_m$ the
results are much more modest.
In 1912 Hardy \cite{hardy_03}, by trying to generalize
Jacobsthal--Vacca's series {\eqref{lk2jd029jde}}
to first generalized Euler's constant, obtained the following series
%
%e39 #&#
\begin{eqnarray}
\label{jh3298uhnjd}
\label{jm9c204dj} \gamma_1 = \frac{\ln2}{2}\sum
_{n=2}^\infty\frac{(-1)^n}{n} \lfloor
\log_2{n}\rfloor\cdot\big(2\log_2{n} - \lfloor
\log_2{2n}\rfloor\big)
\end{eqnarray}
which is, however, not a full generalization of Jacobsthal--Vacca's series
since it also contains irrational coefficients.
In 1924--1927, Kluyver tried, on several occasions \cite{kluyver_03,kluyver_01}, to better Hardy's result and to obtain series for
$\gamma_m$ with
rational terms only, but these attempts were not successful.
There also exist formul{\ae}~similar to Hardy's series.
For instance,
\be\notag
\sum_{n=1}^\infty \frac{\,H^{m}_n - (\gamma+\ln n)^m\,}{n} \,=\,
\begin{cases}
\,-\gamma_1 -\frac{1}{2}\gamma^2+\frac{1}{2}\zeta(2)\,,\qquad & m=1 \\[2mm]
\,-\gamma_2 -\frac{2}{3}\gamma^3 -2\gamma_1\gamma +\frac{5}{3}\zeta(3)\,,\qquad & m=2 \\[2mm]
\,-\gamma_3  -\frac{3}{4}\gamma^4 -3\gamma_2\gamma-3\gamma_1\gamma^2+\frac{43}{8}\zeta(4)\,,\qquad & m=3
\end{cases}
\ee
see e.g.~\cite{furdui_01,mathstack_02}.\footnote{Cases $m=1$
and $m=2$ are discussed in cited references.
Formula for $m=3$ was kindly
communicated to the author by Roberto Tauraso.}Besides, several
asymptotical representations
similar to Hardy's formula
are also known for $\gamma_m$.
For instance, Nikolai M.~G\"unther and Rodion O.~Kuzmin
gave the following formula
%
%e40 #&#
\begin{eqnarray}
\label{d21309dunmd}
\gamma_1\,=\, \sum_{k=1}^n \frac{\,\ln k\,}{k}-\frac{\,\ln^2 \!{n}\,}{2}
-\frac{\,\ln{n}\,}{2\,n} -\frac{\,1-\ln{n}\,}{12\,n^2} + \theta\cdot
\frac{\,11-6\ln{n}\,}{720\,n^4}
\end{eqnarray}
where $0<\theta<1$,
see \cite[\no 388]{gunter_03_eng}.\footnote{In the third edition of
\cite{gunter_03_eng}, published
in 1947, there are two errors in exercise \no388: the value of $\gamma
_1$ is given
as $-0.073927\ldots$ instead of $-0.072815\ldots{}$, and the denominator
of the last term has the value $176$ instead of $720$. These errors
were corrected in the fourth edition of this book, published
in 1951.}
M.~I.~Israilov \cite{israilov_01} generalized expression {\eqref{d21309dunmd}} and showed that the $m$th Stieltjes constant
may be given by a similar semi-convergent asymptotical series
%
%e41 #&#
\be\label{jhx2uxhbcxed}
\gamma_m\,=\,\sum_{k=1}^n \frac{\,\ln^m \! k\,}{k} - \frac{\,\ln
^{m+1} \! n\,}{m+1}
- \frac{\,\ln^m \! n\,}{2n}- \sum_{k=1}^{N-1} \frac{\,{B}_{2k}\,
}{(2k)!}\left[\frac{\ln^m \! x}{x}\right]^{(2k-1)}_{x=n}\!\!
- \theta\cdot\frac{\,{B}_{2N}\,}{(2N)!}\left[\frac{\ln^m\!
x}{x}\right]^{(2N-1)}_{x=n}
\ee
where $m=0, 1, 2,\ldots{}$, $0<\theta<1$, and integers $n$ and $N$
may be arbitrary chosen provided
that $N$ remains finite.\footnote{Note that at fixed
$N$, the greater the number $n$, the more accurate this formula; at $n\to\infty$ it straightforwardly
reduces to {\eqref{k98y9g87fcfcf}}. It seems also appropriate to
note here that, although
G\"unther, Kuzmin and Israilov obtained {\eqref{d21309dunmd}} 
and {\eqref{jhx2uxhbcxed}} independently,  both these formul\ae~may be readily derived from an old 
semi--convergent series for the $\zeta$--function, given, for example, by 
J{\o}rgen P.~Gram in 1895 \cite[p.~304, 2nd unnumbered formula]{gram_01} (this series for
$\zeta(s)$ may be, of course, much older since it is a simple application
of the Euler--Maclaurin summation formula; note also that Gram uses a slightly different convention for the Bernoulli numbers).}\up{,}\footnote{In \cite[Eq.~(3)]{israilov_01},
there is a misprint: in the denominator of the second sum $2r$ should
be replaced by $(2r)!$
[this formula appears correctly on p.~101 \cite{israilov_01}, but with
a misprint in Eq.~(3) on p.~98]. This misprint
was later reproduced in \cite[Theorem~0.3]{eddin_02}.} Using various
series representations
for the $\zeta$-function, it is also possible to obtain corresponding
series for the Stieltjes constants.
For instance, from Ser's series for the $\zeta$-function (see
footnote \ref{kw09h2nb}), it follows that
%
%e42 #&#
\begin{eqnarray}
\gamma_m\,=\,-\frac{1}{\,m+1\,}\sum_{n=0}^\infty\frac{1}{\,n+2\,
}\sum_{k=0}^{n} (-1)^k \binom{n}{k}\frac{\ln^{m+1}(k+1)}{k+1}\,,
\qquad m=0, 1, 2\ldots
\end{eqnarray}
An equivalent result was given by Donal F.~Connon \cite{connon_06}\footnote{Strictly speaking,
Connon \cite[p.~1]{connon_06} gave this formula for the generalized
Stieltjes constants $\gamma_m(v)$,
of~which ordinary Stieltjes constants are simple particular cases
$\gamma_m=\gamma_m(1)$, see e.g.~\cite[Eqs.~(1)--(2)]{iaroslav_07}.}
%
%e43 #&#
\begin{eqnarray}
\gamma_m\,=\,-\frac{1}{\,m+1\,}\sum_{n=0}^\infty\frac{1}{\,n+1\,
}\sum_{k=0}^n (-1)^k \binom{n}{k}\ln^{m+1}(k+1)\,,
\qquad m=0, 1, 2\ldots
\end{eqnarray}
who used Hasse's series for the $\zeta$-function.\footnote{\label{kw09h2nb}The
series representation for $\zeta(s)$,
which is usually attributed to Helmut Hasse, is actually due to Joseph
Ser who derived it in 1926 in a slightly different form.
The series in question is
\be\nonumber
\zeta(s)\,=\,\frac{1}{\,s-1\,}\sum_{n=0}^\infty\frac{1}{\,n+2\,
}\sum_{k=0}^{n} (-1)^k \binom{n}{k}(k+1)^{-s}
\,=\,\frac{1}{\,s-1\,}\sum_{n=0}^\infty\frac{1}{\,n+1\,}\sum
_{k=0}^n (-1)^k \binom{n}{k}(k+1)^{1-s}
\ee
The first variant was given by Ser in 1926 in \cite[p.~1076,
Eq.~(7)]{ser_01}, the second variant was given by Hasse
in 1930 \cite[pp.~460--461]{hasse_01}. The equivalence between two
forms follows from the recurrence relation
for the binomial coefficients.
It is interesting that many writers do not realize that Ser's formula
and Hasse's formula are actually the same
(the problem is also that Ser's paper \cite{ser_01} contains errors,
e.g.~in Eq.~(2), p.~1075, the last term in the second line
should be $(-1)^n (n+1)^{-s}$ instead of $(-1)^n n^{-s}$).
An equivalent series representation was also independently discovered
by Jonathan Sondow
in 1994 \cite{sondow_03}.}
Similarly, using another Ser's series expansions for $\zeta
(s)$,\footnote{Ser's formula \cite[p.~1076, Eq.~(4)]{ser_01},
corrected (see footnote \ref{kw09h2nb}) and written in our
notations, reads
\begin{eqnarray}
\nonumber
\zeta(s)\,=\,\frac{1}{\,s-1\,}+\sum_{n=0}^\infty\big| G_{n+1}\big
|\sum_{k=0}^{n} (-1)^k \binom{n}{k}(k+1)^{-s}.
\end{eqnarray}
} we conclude that
%
%e44 #&#
\begin{eqnarray}
\gamma_m\,=\sum_{n=0}^\infty\big| G_{n+1}\big|\sum_{k=0}^{n}
(-1)^k \binom{n}{k}\frac{\ln^{m}(k+1)}{k+1}\,,
\qquad m=0, 1, 2\ldots
\end{eqnarray}
where $G_n$ are Gregory's coefficients, see
\eqref{njw3uiqch}--\eqref{ldhd9ehn}.\footnote{Note that for $m=0$,
the latter series reduces to Fontana--Mascheroni's series {\eqref{njw3uiqch}}.} This series
was also independently discovered by Marc--Antoine Coppo in 1997 \cite[p.~355, Eq.~(5)]{coppo_01} by the method
of finite differences.
Several more complicated series representations for $\gamma_m$ with
irrational coefficients
may be found in \cite{todd_01,lavrik_01_eng,israilov_01,stankus_01_eng,zhang_01,weisstein_06,coffey_02,coffey_08}.\footnote{Apart 
from formula {\eqref{jhx2uxhbcxed}}, in \cite{israilov_01},
Israilov also obtained
several other series representations for $\gamma_m$.
Stankus \cite{stankus_01_eng} showed that first two Stieltjes
constants may be represented by the series involving the divisor
function. Works of several authors showing that Stieltjes
constants are related to series containing nontrivial zeros of the 
$\zeta$--function are summarized in \cite{weisstein_06}.
Coffey gave various series representations for $\gamma_m$ in \cite{coffey_02,coffey_08}.
However, we also noted that \cite{coffey_02,coffey_08} both 
contain numerous rediscoveries, as well as inaccuracies in attribution of formul\ae~(see also \cite{connon_05}). 
For instance, ``Addison-type series'' \cite[Eq.~(1.3)]{coffey_08} are actually due to Nielsen and Jacobsthal, 
see our Eq.~{\eqref{lk2eojmwjksd}} above.
Numbers $p_{n+1}$ are simply signless Gregory's
coefficients $|G_n|$ and their asymptotics, 
Eq.~(4.8)/(2.82) \cite[p.~473/31]{coffey_02}, is known at
least since 1879 \cite[Eqs.~(25)--(25a)]{schroder_01}
(see also \cite{steffensen_01}, \cite[pp.~106--107]{steffensen_02}). 
Their ``full'' asymptotics, Eq.~(4.10)/(2.84)
\cite[p.~473/31]{coffey_02}, is also known and 
was given, for example, by Van Veen in 1951
\cite{van_veen_01}, \cite{norlund_01}.
Representation (2.17) \cite[p.~455/13]{coffey_02} is
due to Hermite, see Eq.~(13) \cite[p.~541]{iaroslav_07}. 
Formula (1.17) \cite[p.~2052]{coffey_08} was
already obtained at least twice: by Binet in 1839
\cite[p.~257]{binet_01} and by N{\o}rlund in 1923
\cite[p.~244]{norlund_02}, see details in  \cite[Final remark]{iaroslav_08}.
Formula (1.18) \cite[p.~2052]{coffey_08}
straightforwardly follows from Ser's
results \cite{ser_01} dating back to 1926, \emph{etc.}}

%s3.2 #&#
\subsection{Derivation of the series expansion}

\begin{figure}[!t]   
\centering
\includegraphics[width=0.42\textwidth]{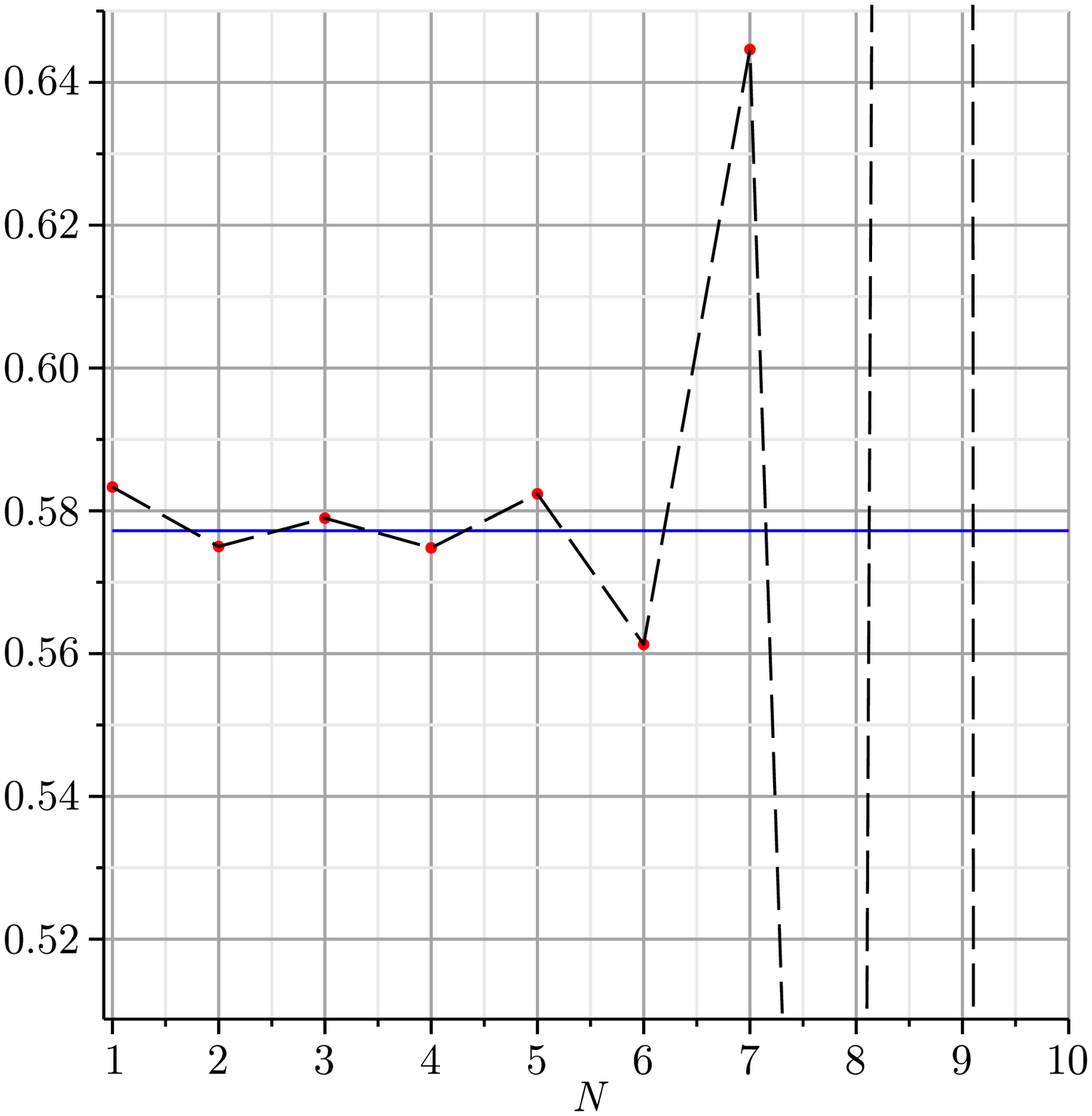}\hspace{8mm}
\includegraphics[width=0.42\textwidth]{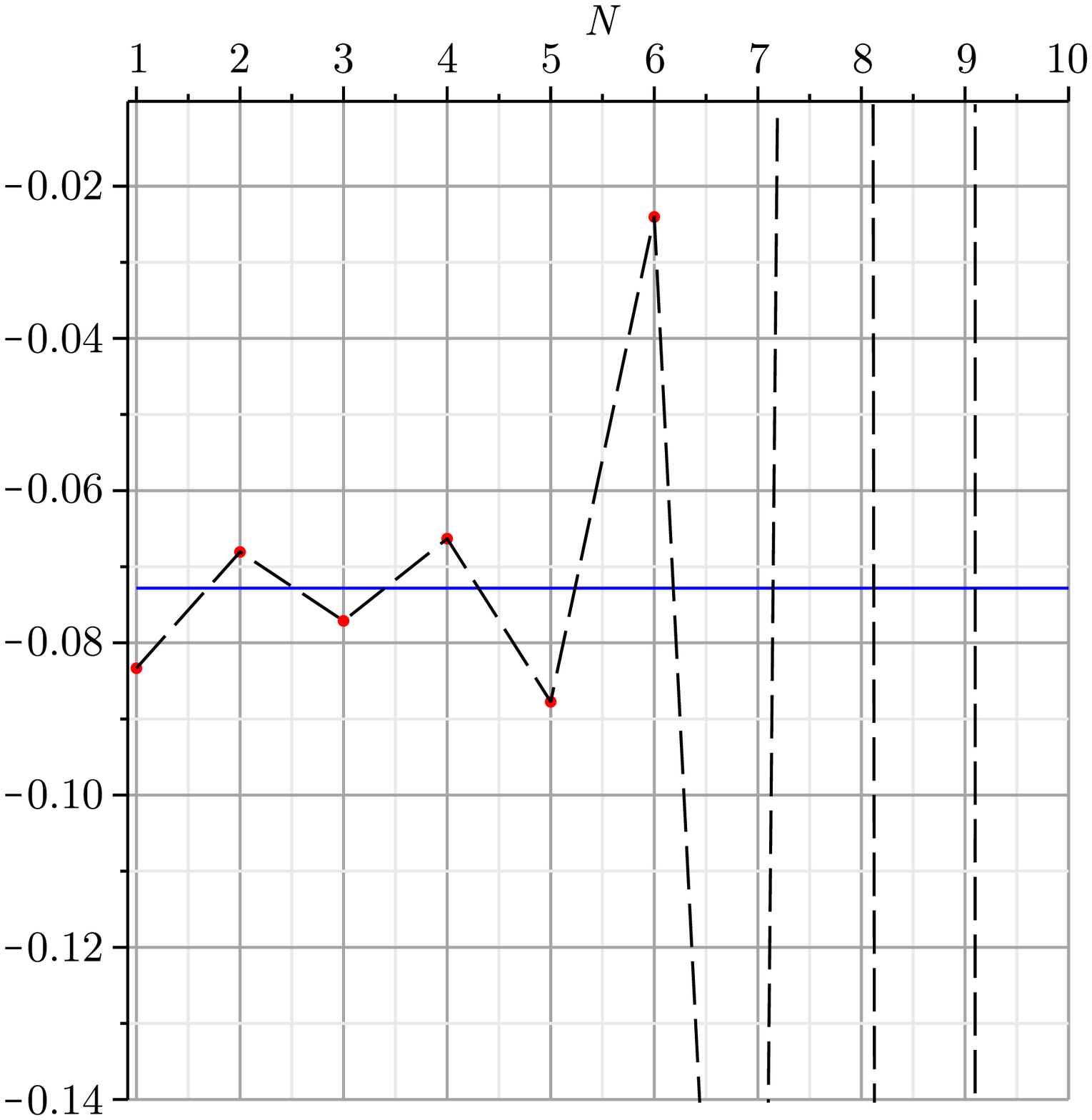}
\caption{Partial sums of series \eqref{c3pf3metdd} for $\gamma$ and $\gamma_1$ (on the left and on the right respectively) 
for $N=1, 2,\ldots, 10$.
Blue lines indicate the true value of $\gamma$ and $\gamma_1$, while dashed black lines with the red points display
corresponding partial sums given by \eqref{c3pf3metdd}.}
\label{ihce239hb}
\end{figure}

Consider now series {\eqref{jkhf3984fhd}}. A formal rearrangement of
this expression
may produce a series with rational terms only for $\gamma_m$.
In view of the fact that
%
%e45 #&#
\begin{eqnarray}
\label{kljd023jdnr}
\zeta(k+1)\,=\,\sum_{n=k}^\infty\frac{\,\big|S_1(n,k)\big|\,
}{n\cdot n!}\,=\,\sum_{n=1}^\infty\frac{\big|S_1(n,k)\big|}{n\cdot
n!}\,,\qquad\quad k=1, 2, 3,\ldots
\end{eqnarray}
see e.g.~\cite[pp.~166, 194--195]{jordan_01},\footnote{See also \cite{shen_01,sato_01},
where this important result was rediscovered much later.
By the way, this formula
may be
generalized to the Hurwitz $\zeta$-function
\begin{eqnarray}
\nonumber
\zeta(k+1,v)\,=\,\sum_{n=k}^\infty\frac{|S_1(n,k)|}{n\cdot(v)_n}\,
,\qquad\quad k=1, 2, 3,\ldots\,,\quad\,
\operatorname{Re}{v}>0\,,
\end{eqnarray}
see \cite{iaroslav_08}. At $n\to\infty$, the
general term of this series is
$O\!\left(\!\dfrac{\,\ln^{k-1}\! n}{n^{v+1}\,}\!\right)$.}
and that
\begin{eqnarray}
\nonumber
\zeta(2k)=(-1)^{k+1}\frac{\,(2\pi)^{2k}\cdot{B}_{2k}\,}{2\cdot
(2k)!}\,=\,
\frac{\,(2\pi)^{2k}\cdot|{B}_{2k}|\,}{2\cdot(2k)!}\,,\qquad\quad
k=1, 2, 3,\ldots
\end{eqnarray}
the formal interchanging of the order of summation in \eqref{jkhf3984fhd} leads to
\be\notag
\begin{array}{l}
\displaystyle
\frac{\,(-1)^m m!\,}{\pi} \sum_{n=1}^\infty\frac{1}{\,n\cdot n!\,} 
\sum_{k=0}^{\lfloor\!\frac{1}{2}n\!\rfloor}\frac{\,(-1)^{k}\big|S_1(2k+2,m+1)\big| \cdot\big|S_1(n,2k+1)\big|}{\,(2\pi)^{2k+1}\,} \asymp\\[8mm]
\displaystyle\qquad\qquad\quad
\,\asymp\,\frac{\,(-1)^m m!\,}{\pi}\sum_{k=0}^\infty\frac{\,(-1)^{k} \big|S_1(2k+2,m+1)\big| \,}{(2\pi)^{2k+1}}
\underbrace{\sum_{n=1}^\infty \frac{\,\big|S_1(n,2k+1)\big|\,}{\,n\cdot n!\,}}_{\zeta(2k+2)} =\\[3mm]
\displaystyle\qquad\qquad\quad
=\,(-1)^m m!\!\sum_{k=1}^{\infty}\frac{\, \big|S_1(2k,m+1)\big| \cdot{B}_{2k}}{(2k)!} 
\end{array}
\ee
Remarking that such operations are not permitted when series are not
absolutely convergent (that is why we wrote $\asymp$
instead of $=$),
we understand why the resulting series diverges. However, since this
series is alternating,
for any prescribed $m$ and $N$, one can always find such $\theta\in
(0,1)$, generally depending
on $m$ and $N$, that
\be\label{349fu3j4nf}
\begin{array}{ll}
\displaystyle
\gamma_m\:&\displaystyle
=\,\frac{1}{2}\delta_{m,0}+(-1)^{m} m!\cdot\!\sum_{k=1}^{N}\frac{\,\big|S_1(2k,m+1)\big|\cdot{B}_{2k}\,}{(2k)!}+ \\[7mm]
&\displaystyle\qquad\qquad\qquad\qquad\qquad\qquad
\,+\, \theta\cdot\frac{\,(-1)^{m} m!\!\cdot \big|S_1(2N+2,m+1)\big|\cdot{B}_{2N+2}\,}{(2N+2)!}=\\[8mm]
&=\begin{cases}
\displaystyle \phantom{+}\frac{1}{2}+ \frac{1}{12}-\frac{1}{120}+\frac{1}{252}-\frac{1}{240}+\frac{1}{132}
% -\frac{691}{32\,760} +\frac{1}{12}
- \ldots \,\quad& m=0 \\[6mm]
\displaystyle -\frac{1}{12}+\frac{11}{720}-\frac{137}{15\,120}+\frac{121}{11\,200}-\frac{7\,129}{332\,640}+\frac{57\,844\,301}{908\,107\,200}
% -\frac{1\,145\,993}{4\,324\,320}
-\ldots \,\quad& m=1 \\[6mm]
\displaystyle +0-\frac{1}{60}+\frac{5}{336}-\frac{469}{21\,600}+\frac{6\,515}{133\,056}-\frac{131\,672\,123}{825\,552\,000}
 +\frac{63\,427}{89\,100}
-\ldots \,\quad& m=2 \\[6mm]
\displaystyle -0+\frac{1}{120}-\frac{17}{1\,008}+\frac{967}{28\,800}-\frac{4\,523}{49\,896}+\frac{33\,735\,311}{101\,088\,000}
-\frac{9\,301\,169}{5\,702\,400} 
+ \ldots \,\quad& m=3 
\end{cases}
\end{array}
\ee
holds strictly.\footnote{There is another way to
obtain {\eqref{349fu3j4nf}}. Consider
first {\eqref{jhx2uxhbcxed}} at $n=1$, 
and then use {\eqref{iu2d092n1}}
to show that \break $\,\left[\frac{d^{n}}{dx^n}\frac{\ln^m \! x}{x}\right]_{x=1}\!=m!\, S_1(n+1,m+1)\,$.}
Moreover, taking into account {\eqref{uf87tfuy89}}, above series may be
always written in a form without Stirling numbers.
For instance, for Euler's constant and for first three Stieltjes
constants, it becomes
\be\label{c3pf3metdd}
\begin{array}{ll}
\displaystyle
\gamma\,
=\,+\frac{1}{2}\left\{1+\sum_{k=1}^{N}\frac{\,{B}_{2k}}{\,k\,} 
+\, \theta\cdot\frac{\,{B}_{2N+2}\,}{\,N+1\,}\right\}\\[7mm]
\displaystyle
\gamma_1\,
=\,-\frac{1}{2}\sum_{k=1}^{N}\frac{\,{B}_{2k}\cdot H_{2k-1}\,}{\,k\,} 
+\, \theta\cdot\frac{\,{B}_{2N+2}\cdot H_{2N+1}\,\,}{\,2N+2\,}
\\[7mm]
\displaystyle
\gamma_2\,
=\,+\frac{1}{2}\sum_{k=1}^{N}\frac{\,{B}_{2k}\cdot \big\{H^2_{2k-1} - H^{(2)}_{2k-1}\big\}\,}{\,k\,} 
+\, \theta\cdot\frac{\,{B}_{2N+2}\cdot \big\{H^2_{2N+1} - H^{(2)}_{2N+1}\big\}\,\,}{\,2N+2\,}
\\[7mm]
\displaystyle
\gamma_3\,
=\,-\frac{1}{2}\sum_{k=1}^{N}\frac{\,{B}_{2k}\cdot \big\{H^3_{2k-1} - 3H_{2k-1} H^{(2)}_{2k-1}+2H^{(3)}_{2k-1}\big\}\,}{\,k\,}\,+ \\[4mm]
\displaystyle\qquad\qquad\qquad\qquad\quad
+\, \theta\cdot\frac{\,{B}_{2N+2}\cdot \big\{H^3_{2N+1} - 3H_{2N+1} H^{(2)}_{2N+1}+2H^{(3)}_{2N+1}\big\}\,\,}{\,2N+2\,}
\end{array}
\ee
where $0<\theta<1$ and $ N<\infty$ (these parameters are different in
all equations, and in each equation $\theta$, in general, depends on $N$).
By the way, one may notice that first of these formul{\ae}~coincide
with Stirling series {\eqref{lkce02m}}, while other formul{\ae}~are, to
our knowledge, new and
seem to be never released before.
Derived formal series are alternating and are sometimes referred to as
\emph{semi-convergent series} or \emph{divergent enveloping series}.\footnote{These series were 
an object of study of almost all great
mathematicians; the reader interested
in a deeper study of these series may wish
to consult the following literature: \cite{borel_01,erdelyi_01,hardy_02}, \cite[Chapt.~XI]{bromwich_01},
\cite[Chapt.~4, \S1]{polya_01_eng},
\cite{knopp_01},
exercises \no374--\no388 in \cite[pp.~46--48]{gunter_03_eng}, \cite{whittaker_01,paplauskas_01_eng,vorobiev_01,copson_01,dingle_01,olver_01,bleistein_01,ramis_01,ramis_02,euler_03,malgrange_01}.}
It is also easy to see that they diverge very rapidly
\begin{eqnarray*}
\displaystyle
\frac{\,\big|S_1(2k,m+1)\big|\cdot{B}_{2k}}{(2k)!\,}\,&\sim&\,
2(-1)^{k-1}\frac{\,(2k-1)!\cdot\ln^m(2k-1)\,}{m!\cdot(2\pi)^{2k}}
\\[5mm]
\,&\sim&\,\frac{\,2\sqrt{\pi\,}(-1)^{k-1}}{m!}\cdot\frac{\,\ln^m k\,
}{\,\sqrt{k\,}\,}\cdot \left(\frac{k}{\pi e}\right)^{2k}\,,
\end{eqnarray*}
$k\to\infty$, $m=0, 1, 2,\ldots{}$,
so rapidly that even the corresponding
power series $\sum\frac{\,|S_1(2k,m+1)|\,}{(2k)!}{B}_{2k}x^{2k}$
diverges everywhere.\footnote{Coefficients $\big|S_1(2k,m+1)\big|$
and Bernoulli numbers both
grow very quickly: as $k\to\infty$ we have $\big|S_1(2k,m+\nobreak 1)\big
|\sim(2k-1)!\ln^m(2k-1)/m!$, see {\eqref{j6s8g64r}}, and
${B}_{2k}\sim2(-1)^{k-1}(2\pi)^{-2k}(2k)!$, see e.g.~\cite[p.~5]{krylov_01}, \cite[p.~261]{gelfond_01}.}
Behaviour of this series for first two Stieltjes constants is shown in
{Fig.~\ref{ihce239hb}}.
%%%%%%%%%%%%%%%%%%%%%%%%%%%%%%%%%%%%%%%%%%%%%%%%%%%%%%%%
%%%%%%%%%%%%%%%%%%%%%%%%%%%%%%%%%%%%%%%%%%%%%%%%%%%%%%%

%s3.3 #&#
\subsection{Some series transformations applied to the derived
divergent series}
\begin{figure}[!t]   
\centering
\includegraphics[width=0.42\textwidth]{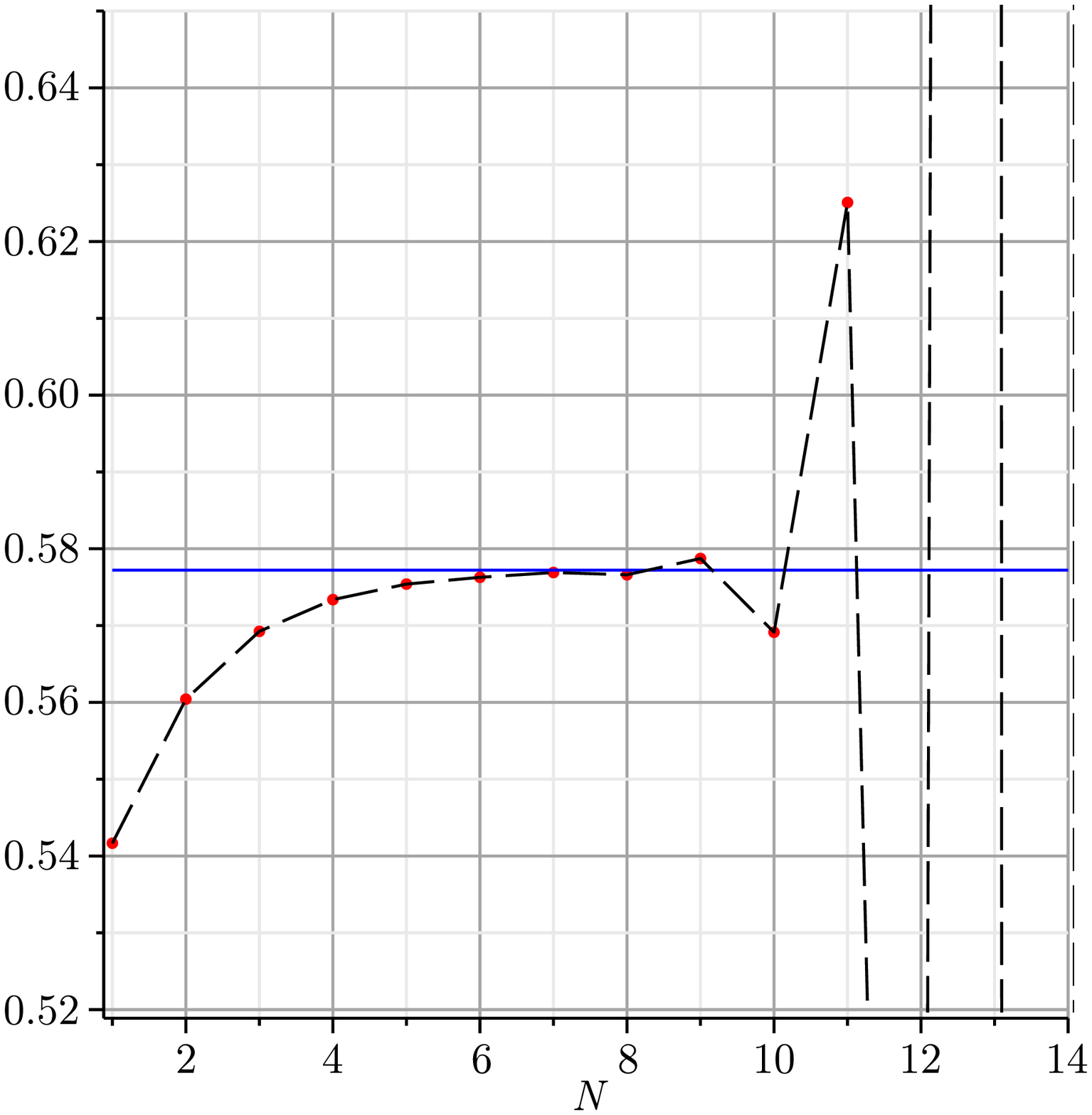}\hspace{8mm}
\includegraphics[width=0.42\textwidth]{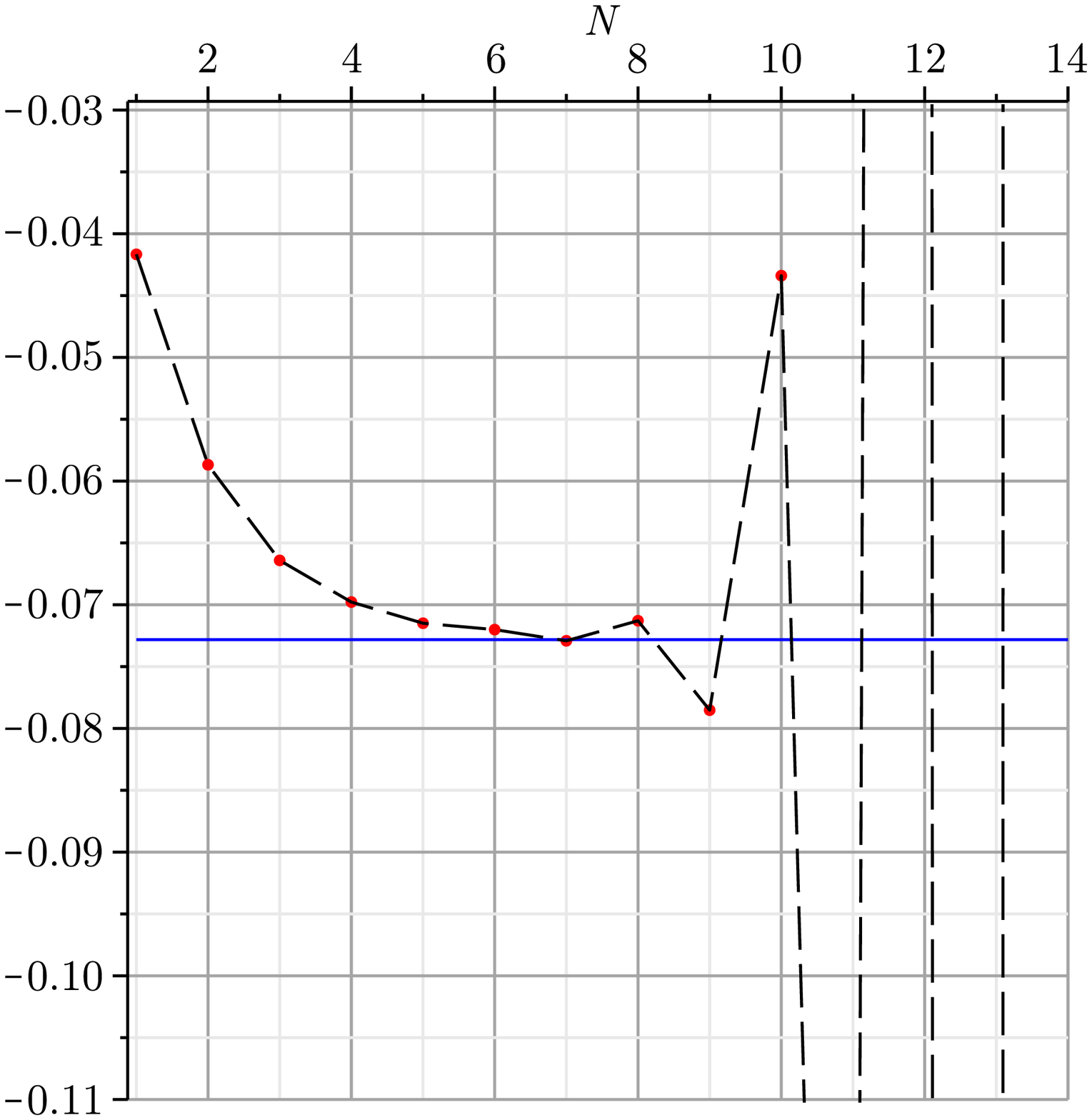}
\caption{Euler's transformation of series \eqref{c3pf3metdd} for $\gamma$ and $\gamma_1$ (on the 
left and on the right respectively) for $N=1, 2,\ldots, 14$.
Blue lines indicate the true value of $\gamma$ and $\gamma_1$, while dashed black lines with the red points display
corresponding partial sums given by \eqref{kjwhe932hs}.}
\label{lkf3rh0hjda}
\end{figure}

In order to convert {\eqref{349fu3j4nf}} into a convergent series, one
may try to apply various series transformations and regularization procedures.
However, since {\eqref{349fu3j4nf}} is strongly divergent, the use of
standard summation methods may not result in a convergent series. For example,
applying Euler's transformation\footnote{See e.g.~\cite[pp.~244--246]{knopp_01}, \cite[pp.~144 \& 170--171]{knopp_02},
\cite[pp.~269--278 \& 305--306]{vorobiev_01}.} we obtain another
series with rational coefficients only\\[4mm]
\be\notag
\begin{array}{ll}
\displaystyle
\gamma_m\, =\,\frac{1}{2}\delta_{m,0}+(-1)^{m} m!\!\sum_{k=1}^N \frac{1}{\,2^k}\sum_{n=1}^{k} 
\frac{\,\big|S_1(2n,m+1)\big|\cdot{B}_{2n}}{\,(2n)!\,}\cdot\!\binom{k-1}{n-1} + R_m(N)=\\[8mm]
\end{array}
\ee

\be\label{kjwhe932hs}
\begin{array}{ll}
\displaystyle
=\begin{cases}
\displaystyle \phantom{+} \frac{1}{2}+\frac{1}{24}+\frac{3}{160}+\frac{89}{10\,080}+\frac{37}{8\,960}
+\frac{299}{147\,840}
% +\frac{323}{366\,080}
+ \ldots\,\quad& m=0 \\[5mm]
\displaystyle -\frac{1}{24}-\frac{49}{2\,880}-\frac{187}{24\,192}-\frac{5\,431}{1\,612\,800}-\frac{91\,151}{53\,222\,400}
% -\frac{29\,947\,039}{58\,118\,860\,800}
-\ldots \,\quad& m=1 \\[5mm]
\displaystyle -\frac{1}{240}-\frac{31}{13\,440}-\frac{4\,093}{2\,419\,200}-\frac{50\,789}{106\,444\,800}-\frac{602325403}{581188608000}
% +\frac{36\,800\,663}{26\,417\,664\,000}
+\ldots \,\quad& m=2 \\[5mm]
\displaystyle +\frac{1}{480}-\frac{1}{40\,320}+\frac{1\,609}{3\,225\,600}-\frac{120\,749}{159\,667\,200}+\frac{694\,773\,847}{498\,161\,664\,000}
% -\frac{1\,349\,033\,801}{332\,107\,776\,000}
- \ldots \,\quad& m=3 
\end{cases} 
\vspace{2mm}
\end{array}
\ee
which are all divergent (their remainder $R(N)\to\infty$ as $N\to
\infty$).
At the same time, these series behave much better than {\eqref{349fu3j4nf}}. Thus,
the series for $\gamma$ starts to clearly diverge only from $N\geqslant
10$, and that
for $\gamma_1$ only from $N\geqslant8$, see {Fig.~\ref{lkf3rh0hjda}}.
The minimum error for the first series corresponds to $N=7$ and equals
$3\times10^{-4}$,
that for $\gamma_1$ also corresponds to $N=7$ and equals $9\times10^{-5}$.
Attempts to regularize series \eqref{349fu3j4nf}--\eqref{c3pf3metdd}
with the help of Ces\`aro summation are also fruitless
since its general term grows more rapidly than $k$ at $k\to\infty$.
Similarly, Borel's summation does not provide
a convergent result.
%%%%%%%%%%%%%%%%%%%%%%%%%%%%%%%%%%%%%%%%%%%%%%%%%%%%%%%%%%
%%%%%%%%%%%%%%%%%%%%%%%%%%%%%%%%%%%%%%%%%%%%%%%%%%%%%%%%%

%s3.4 #&#
\subsection{An estimate for generalized Euler's constants}
Finally, we note that derived formal series {\eqref{349fu3j4nf}}
provides an estimation
for generalized Euler's constants. Since this series is enveloping, its
true value always lies between
two neighbouring partial sums. If, for example, we retain only
the first non-vanishing term, the $m$th Stieltjes constant
will be stretched between the first and second non-vanishing partial sums.
Thus, accounting for known properties of the Stirling numbers of the
first kind
\be\notag
\begin{array}{ll}
\big|S_1(m+1,m+1)\big|\,=\,1 \\[3mm]
\big|S_1(m+2,m+1)\big|\,=\,\frac{1}{2}(m+1)(m+2) \\[3mm]
\big|S_1(m+3,m+1)\big|\,=\,\frac{1}{24}(m+1)(m+2) (m+3)(3m+8) \\[-3mm]
\end{array}
\ee

\be\notag
\big|S_1(m+4,m+1)\big|\,=\,\tfrac{1}{48}(m+1)(m+2) (m+3)^2(m+4)^2
\ee
which are valid for $m=1, 2, \ldots\,$, we have
\begin{equation}\label{lkx2xmx}
\begin{array}{ll}
\displaystyle-\frac{\,\big|{B}_{m+1}\big|\,}{m+1} < \gamma_m  <
\frac{\,(3m+8)\cdot\big|{B}_{m+3}\big|\,}{24} - \frac{\,\big|{B}_{m+1}\big|\,}{m+1}  , 
& m=1, 5, 9,\ldots\\[18pt]
\displaystyle 
\frac{\,\big|{B}_{m+1}\big|\,}{m+1} - \frac{\,(3m+8)\cdot\big|{B}_{m+3}\big|\,}{24}
< \gamma_m  < \frac{\,\big|{B}_{m+1}\big|\,}{m+1} , & m=3, 7, 11,\ldots\\[18pt]
\displaystyle -\frac{\,\big|{B}_{m+2}\big|\,}{2} < \gamma_m
  <  \frac{\,(m+3)(m+4)\cdot\big|{B}_{m+4}\big|\,}{48} - \frac{\,\big|{B}_{m+2}\big|\,}{2} ,
 \qquad &  m=2, 6, 10, \ldots\\[18pt]
\displaystyle 
\frac{\,\big|{B}_{m+2}\big|\,
}{2} - \frac{\,(m+3)(m+4)\cdot\big|{B}_{m+4}\big|\,}{48}
< \gamma_m   < \frac{\,\big|{B}_{m+2}\big|\,}{2}, & m=4, 8, 12, \ldots\\
\end{array}
\end{equation}
Case $m=4, 8, 12, \ldots$ may be also extended to $m=0$, if we recall that in this case it gives bounds 
for $\gamma-\frac12$, see {\eqref{349fu3j4nf}}. This case yields 
$\,\frac{23}{40}<\gamma<\frac{7}{12}\,$, which is undoubtedly true. Note also that above
bounds are always rational, which may be of interest in certain circumstances.

Estimation \eqref{lkx2xmx} is relatively tight for moderate values of $m$, and
becomes less and less accurate as $m$ increases.
However, even for large $m$, it remains more accurate than the
well-known Berndt's estimation
\be
\big|\gamma_m\big|\,\leqslant\,
\begin{cases}
\displaystyle \frac{2\,(m-1)!}{\pi^m}\,,\qquad &  m=1, 3, 5,\ldots \\[4mm]
\displaystyle \frac{4\,(m-1)!}{\pi^m}\,,\qquad &  m=2, 4, 6,\ldots 
\end{cases}
\ee
see \cite[pp.~152--153]{berndt_02}, more accurate than Lavrik's
estimation $|\gamma_m|\leqslant m! \, 2^{-m-1}$, see
\cite[Lemma~4]{lavrik_01_eng}, more accurate than Israilov's estimation
\be
|\gamma_m|\leqslant \,\frac{ m!\, C(k)\,}{(2k)^{m}}
\ee
for $k=1,2,3$, where $C(1)=\frac
{1}{2}$, $C(2)=\frac{7}{12}$, $C(3)=\frac{11}{3}$,
see \cite{israilov_01,adell_01}, and more accurate than
Nan-You--Williams' estimation
\be
\big|\gamma_m\big|\,\leqslant\,
\begin{cases}
\displaystyle \frac{2\,(2m)!}{m^{m+1}(2\pi)^m}\,,\qquad &  m=1, 3, 5,\ldots \\[4mm]
\displaystyle \frac{4\,(2m)!}{m^{m+1}(2\pi)^m}\,,\qquad &  m=2, 4, 6,\ldots 
\end{cases}
\ee
see \cite[pp.~148--149]{zhang_01}.
Besides, our estimation also contains a sign, while above estimations
are signless.
At the same time, {\eqref{lkx2xmx}} is worse than
Matsuoka's estimation \cite{matsuoka_01,matsuoka_02}, $
|\gamma_m|<10^{-4}\ln^m  m$, $m\geqslant5$,
which, as far as we know, is currently the best known estimation in
terms of elementary functions
for the Stieltjes constants.\footnote{Numerical simulations suggest
that Matsuoka's estimation \cite{matsuoka_01,matsuoka_02} may
be considerably improved, see e.g.~\cite{kreminski_01}. Recently,
Knessl and Coffey reported
that they succeeded to significantly better Matsuoka's estimation and
even to predict the sign of $\gamma_m$. The authors published
their findings in \cite{knessl_01}, and also reprinted them in \cite[Theorem~1]{knessl_02}. We, however, were not able to verify these results,
because several important details related to $v(n)$ from pp.~179--180
\cite{knessl_01} were omitted.
Estimation of the similar nature was later proposed by Adell \cite{adell_01}, but Saad Eddin \cite[Tab.~2]{eddin_01}, \cite{eddin_02} reported
that Adell's estimation may provide less accurate results than
Matsuoka's estimation. Saad Eddin also provides an interesting
estimation for the Stieltjes constants, see \cite[Tab.~2]{eddin_01}
and mentions some further works related to the estimations
of the derivatives of certain $L$-functions. Yet, very recently we
found another work devoted to
the estimation of Stieltjes constants \cite{ahmed_01}; the latter
resorts to the Lambert $W$-function.\\[-17mm]}
Note, however, that estimation's bounds {\eqref{lkx2xmx}}
may be bettered if we transform parent series {\eqref{349fu3j4nf}}
into a less divergent series, provided the new series remains enveloping.

\section*{Acknowledgments} 
The author is grateful to Yaroslav Naprienko for his remarks and comments.

\vspace{165mm}

% \konecpar

%%%%%%%%%%%%%%%%%%%%%%%%%%%%%%%%%%%%%%%%%%%%%%%%%%%%%%%%%%%%
\begin{figure}[!h]   
\centering
\includegraphics[width=0.3\textwidth]{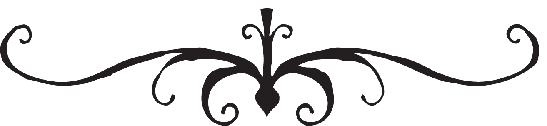}
\end{figure}
%%%%%%%%%%%%%%%%%%%%%%%%%%%%%%%%%%%%%%%%%%%%%%%%%%%%%%%%%%%

\vspace{35mm}

%%%%%%%%%%%%%%%%%%%%%%%%%%%%%%%%%%%%%%%%%%%%%%%%%%%%%%%%%%%%%%%%%%%%%%%%%%%%%%%%%%%%%%%%%%%%%%%%%%%%%%%%%%%%%%%%%
%%%%%%%%%%%%%%%%%%%%%%%%%%%%%%%%%%%%%%%%%%%%%%%%%%%%%%%%%%%%%%%%%%%%%%%%%%%%%%%%%%%%%%%%%%%%%%%%%%%%%%%%%%%%%%%%%
%%%%%%%%%%%%%%%%%%%%%%%%%%%%%%%%%%%%%%%%%%%%%%%%%%%%%%%%%%%%%%%%%%%%%%%%%%%%%%%%%%%%%%%%%%%%%%%%%%%%%%%%%%%%%%%%%
%%%%%%%%%%%%%%%%%%%%%%%%%%%%%%%%%%%%%%%%%%%%%%%%%%%%%%%%%%%%%%%%%%%%%%%%%%%%%%%%%%%%%%%%%%%%%%%%%%%%%%%%%%%%%%%%%
\appendix
\section{Two simple integral formul\ae~for the Stirling numbers of the first kind}\label{ijd20jdk2nd}
\subsection{The first formula}
Consider equation (\ref{x2l3dkkk03d}a) defining signless Stirling
numbers of the first kind.
Dividing both sides by $z^{k+1}$, where $k=1, 2, 3,\ldots{}$, and
integrating along a
simple closed curve $L$ encircling the origin in the counterclockwise
direction, we have
\be
\ointctrclockwise\limits_{L}\!\!\frac{(z)_n}{z^{k+1}}\, dz\,=
\sum_{l=1}^n \big|S_1(n,l) \big|\cdot \underbrace{\ointctrclockwise\limits_{L} z^{l-k-1}\, dz}_{2\pi i \, \delta_{l,k}}
=\,2\pi i \, \big|S_1(n,k) \big|
\ee
in virtue of Cauchy's theorem.
In practice, it is common to take as $L$ the unit circle, and hence
\be
\big|S_1(n,k) \big|\,=\,\frac{1}{\,2\pi i\,}\!\!\ointctrclockwise\limits_{|z|=1}\!\!\frac{(z)_n}{z^{k+1}}\, dz\,
=\,\frac{1}{\,2\pi\,}\!\!\! \int\limits_{\alpha-\pi}^{\alpha+\pi}\!\! \frac{\,\Gamma
(n+e^{i\varphi})\,}{\Gamma(e^{i\varphi})}\, e^{-i\varphi k}\,
d\varphi
\ee
where $\alpha\in\mathbb{R}$.

Alternatively, in virtue of the Cauchy residue theorem, we also have
%
%eA.3 #&#
\be
\big|S_1(n,k) \big|\,=\,\res_{z=0}\frac{\Gamma(n+z)}{\,z^{k+1}\Gamma(z)\,}\,,\qquad\quad n,k\,=\,1,2,3,\ldots
\ee

%sA.2 #&#
\subsection{The second formula}
Consider generating equation (\ref{ld2jr3mnfdmd}a) for the unsigned
Stirling numbers of the first kind.
Proceeding as above and then making a change of variable
$z=re^{i\varphi}$, we obtain
\be\label{ock2w3jkmd1}
\big|S_1(n,k) \big|\, =\,\frac{(-1)^k}{\,2\pi i\,}\cdot\frac{n!}{k!}\!\!\ointctrclockwise\limits_{|z|=r}\!\!\frac{\ln^{k}(1-z)}{z^{n+1}}\, dz\,
=\,\frac{(-1)^k}{\,2\pi\,}\cdot\frac{n!}{k!}\!\!\int\limits_{\alpha-\pi}^{\alpha+\pi}\!\!
\frac{\,\ln^k\!\big(1-re^{i\varphi}\big)\,}{r^n} \, e^{-i\varphi n}\, d\varphi\, 
\ee
where, due to the radius of convergence of (\ref{ld2jr3mnfdmd}a),
$0<r<1$, and $\alpha\in\mathbb{R}$.
Obviously,
the line integral in the middle may be taken not only along the
indicated circle, but along any simple closed curve
encircling the origin in the counterclockwise direction and lying inside
the unit circle $|z|=1$.
By the way, it is interesting that
the same integral taken between $(0,1)$ reduces to a finite linear combination
of Stirling numbers of the first kind and $\zeta$-functions, see
\cite[Sect.~2.2, Eq.~(48)]{iaroslav_08}.
It seems also appropriate to note here that several slightly different
integral formul{\ae}~of the same kind as {\eqref{ock2w3jkmd1}}
were given by Dingle \cite[pp.~92, 199]{dingle_01}. The author, however, did
not specify the integration contour. This inaccuracy
has been partially corrected by Temme \cite[p.~237]{temme_02},
who indicated that the integration contour should be a ``small
circle''.\footnote{We, however, note this condition is not really necessary;
it is sufficient that the integration path be a simple closed curve
encircling the origin in the right direction and
lying inside the unit disc. This remark may be important for the numerical
evaluation of these integrals, which have been reported as not
well-suited for these purposes.}

Similarly to the previous case, the signless Stirling numbers of the
first kind may be also given by the following residue
\be
\big|S_1(n,k) \big|\,= \,(-1)^k \frac{n!}{k!}\res_{z=0}\frac{\ln^{k}(1-z)}{z^{n+1}}\,,\qquad\quad n,k\,=\,1,2,3,\ldots
\ee

\begin{figure}[!t]   
\centering
\includegraphics[width=0.7\textwidth]{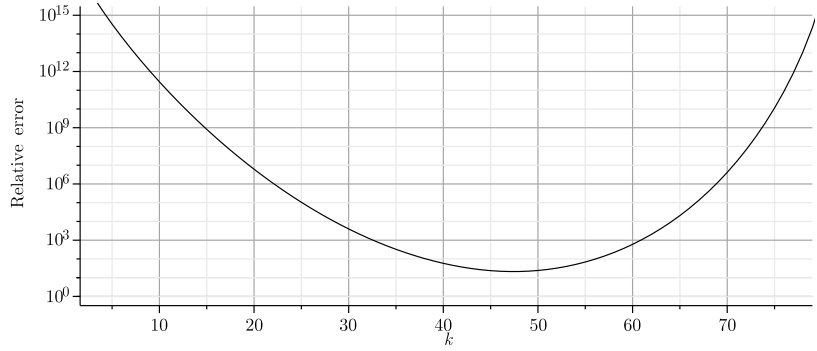}\vspace{-6mm}
\caption{Relative error of the estimation for the Stirling numbers of the first kind given by \protect\eqref{ock2w3jkmd}
as a function of $k$ for $n=80$, logarithmic scale.}
\label{cvkw0er9j}
\end{figure}

Finally, expression {\eqref{ock2w3jkmd1}} may be also used to readily
get an interesting bound for Stirling numbers.
In fact, one may notice that if $r\leqslant1-e^{-1}\approx0.63$, then the
principal branch $\big|\ln^{k}(1-z)\big|\leqslant1$
independently of $k$. Therefore, for $k=1,2,\ldots,n$
\be\label{ock2w3jkmd}
\big|S_1(n,k) \big|\,\leqslant\,\frac{n!}{\,2\pi\,k!}\!\!\ointctrclockwise\limits_{|z|=r}
\!\!\left|\frac{\ln^{k}(1-z)}{z^{n+1}}\right|\cdot \big|dz\big|\, 
\leqslant\,\frac{n!}{\,\big(1-e^{-1}\big)^n\,k!\,}
\ee
since $\big|dz\big|=r\,d\varphi$. If $n$ is prescribed, this
estimation is relatively rough for small and large
(close to $n$) factors $k$; in contrast, for values of $k$ which are
slightly greater than~$n/2$,
this estimation is quite accurate,
see {Fig.~\ref{cvkw0er9j}}.
%%%%%%%%%%%%%%%%%%%%%%%%%%%%%%%%%%%%%%%%%%%%%%%%%%%%%%%%
%%%%%%%%%%%%%%%%%%%%%%%%%%%%%%%%%%%%%%%%%%%%%%%%%%%%%%%


\begin{thebibliography}{100}
\providecommand{\url}[1]{\texttt{#1}}
\providecommand{\urlprefix}{URL }
\expandafter\ifx\csname urlstyle\endcsname\relax
  \providecommand{\doi}[1]{doi:\discretionary{}{}{}#1}\else
  \providecommand{\doi}{doi:\discretionary{}{}{}\begingroup
  \urlstyle{rm}\Url}\fi

\bibitem{stieltjes_01}
Correspondance d'Hermite et de Stieltjes. Vol.~1 and 2, Gauthier-Villars,
  Paris, 1905.

\bibitem{ramanujan_01}
Collected papers of Srinivasa Ramanujan, Cambridge, 1927.

\bibitem{abramowitz_01}
\textit{M.~Abramowitz} and \textit{I.~A. Stegun}, Handbook of mathematical
  functions with formula, graphs and mathematical tables [Applied mathematics
  series \no 55], US Department of Commerce, National Bureau of Standards,
  1961.

\bibitem{adamchik_03}
\textit{V.~Adamchik}, On \text{Stirling} numbers and \text{Euler's} sums,
  Journal of Computational and Applied Mathematics, vol.~79, pp.~119--130
  (1997).

\bibitem{addison_01}
\textit{A.~W. Addison}, A series representation for \text{Euler's} constant,
  The American Mathematical Monthly, vol.~74, pp.~823--824  (1967).

\bibitem{adelberg_01}
\textit{A.~Adelberg}, 2-\text{Adic} congruences of \text{Nörlund} numbers and
  of \text{Bernoulli} numbers of the second kind, Journal of Number Theory,
  vol.~73, no.~1, pp.~47--58  (1998).

\bibitem{adell_01}
\textit{J.~A. Adell}, Asymptotic estimates for \text{Stieltjes} constants: a
  probabilistic approach, Proceedings of the Royal Society A, vol.~467,
  pp.~954--963  (2011).

\bibitem{ainsworth_01}
\textit{O.~R. Ainsworth} and \textit{L.~W. Howell}, An integral representation
  of the generalized \text{Euler--Mascheroni} constants, NASA Technical paper
  2456, pp.~1--11  (1985).

\bibitem{alabdulmohsin_01}
\textit{I.~M. Alabdulmohsin}, Summability calculus, arXiv:1209.5739v1  (2012).

\bibitem{appel_01}
\textit{P.~Appel}, Développement en série entière de $(1+ax)^{\frac{1}{x}}$,
  Archiv der Mathematik und Physik, vol.~65, pp.~171--175  (1880).

\bibitem{arakawa_01}
\textit{T.~Arakawa}, \textit{T.~Ibukiyama} and \textit{M.~Kaneko}, Bernoulli
  Numbers and Zeta Functions, Springer Monographs in Mathematics, Japan, 2014.

\bibitem{barrow_01}
\textit{D.~F. Barrow}, Solution to problem \no 4353, American Mathematical
  Monthly, vol.~58, pp.~116--117  (1951).

\bibitem{bateman_01}
\textit{H.~Bateman} and \textit{A.~Erdélyi}, Higher Transcendental Functions
  [in 3 volumes], Mc Graw--Hill Book Company, 1955.

\bibitem{bellavista_01}
\textit{L.~V. Bellavista}, On the \text{Stirling} numbers of the first kind
  arising from probabilistic and statistical problems, Rendiconti del Circolo
  Matematico di Palermo, vol.~32, no.~1, pp.~19--26  (1983).

\bibitem{bender_01}
\textit{E.~A. Bender} and \textit{S.~G. Williamson}, Foundations of
  Combinatorics with Applications, Addison--Wesley, USA, 1991.

\bibitem{berndt_02}
\textit{B.~C. Berndt}, On the \text{Hurwitz Zeta--function}, Rocky Mountain
  Journal of Mathematics, vol.~2, no.~1, pp.~151--157  (1972).

\bibitem{binet_01}
\textit{M.~J. Binet}, Mémoire sur les intégrales définies eulériennes et
  sur leur application à la théorie des suites, ainsi qu'à l'évaluation des
  fonctions des grands nombres, Journal de l'École Royale Polytechnique, tome
  XVI, cahier 27, pp.~123--343  (1839).

\bibitem{iaroslav_08}
\textit{I.~V. Blagouchine}, Two series expansions for the logarithm of the
  gamma function involving \text{Stirling} numbers and containing only rational
  coefficients for certain arguments related to $\pi^{-1}$, Journal of Mathematical 
  Analysis and Applications (Elsevier), vol.~442, no.~2, pp.~404--434 (2016), arXiv:1408.3902 (2014).

\bibitem{iaroslav_07}
\textit{I.~V. Blagouchine}, A theorem for the closed--form evaluation of the
  first generalized \text{Stieltjes} constant at rational arguments and some
  related summations, Journal of Number Theory (Elsevier), vol.~148,
  pp.~537--592 and vol.~151, pp.~276--277  (2015), arXiv:1401.3724 (2014).

\bibitem{bleistein_01}
\textit{N.~Bleistein} and \textit{R.~A. Handelsman}, Asymptotic Expansions of
  integrals, Holt, Rinehart and Winston, USA, 1975.

\bibitem{boole_01}
\textit{G.~Boole}, Calculus of finite differences (edited by J.~F.~Moulton, 4th
  ed.), Chelsea Publishing Company, New-York, USA, 1957.

\bibitem{borel_01}
\textit{É.~Borel}, Leçons sur les séries divergentes (2nd Edn.),
  Gauthier--Villars, Paris, France, 1928.

\bibitem{briggs_01}
\textit{W.~E. Briggs}, The irrationality of $\gamma$ or of sets of similar
  constants, Vid. Selsk. Forh. (Trondheim), vol.~34, pp.~25--28  (1961).
  
\bibitem{briggs_02}
\textit{W.~E. Briggs} and \textit{S. Chowla}, The power series coefficients of $\zeta(s)$, 
American Mathematical Monthly, vol.~62, no.~5, pp.~323--325  (1955).

\bibitem{bromwich_01}
\textit{T.~J. \text{I'a} Bromwich}, An introduction to the theory of infinite
  series, Macmillan and Co. Limited, St--Martin Street, London, 1908.

\bibitem{butzer_01}
\textit{P.~M. Butzer} and \textit{M.~Hauss}, \text{Stirling} functions of the
  first and second kinds; some new applications, In Israel Mathematical
  Conference Proceedings, Approximation, Interpolation, and Summability, in
  Honor of Amnon Jakimovski on his Sixty-Fifth Birthday (Ed. S. Baron and D.
  Leviatan), Bar-Ilan University, Tel Aviv, June 4--8 1990, pp.~89--108,
  Weizmann Press, Israel  (1991).

\bibitem{butzer_02}
\textit{P.~M. Butzer}, \textit{C.~Markett} and \textit{M.~Schmidt},
  \text{Stirling} numbers, central factorial numbers, and representations of
  the Riemann zeta function, Results in mathematics, vol.~19, no.~3--4,
  pp.~257--274  (1991).

\bibitem{candelpergher_01}
\textit{B.~Candelpergher} and \textit{M.-A. Coppo}, A new class of identities
  involving \text{Cauchy} numbers, harmonic numbers and zeta values, The
  Ramanujan Journal, vol.~27, pp.~305--328  (2012).

\bibitem{carlitz_02}
\textit{L.~Carlitz}, Some theorems on \text{Bernoulli} numbers of higher order,
  Pacific Journal of Mathematics, vol.~2, no.~2, pp.~127--139  (1952).

\bibitem{carlitz_01}
\textit{L.~Carlitz}, A note on \text{Bernoulli} and \text{Euler} polynomials of
  the second kind, Scripta Mathematica, vol.~25, pp.~323--330  (1961).

\bibitem{carlitz_03}
\textit{L.~Carlitz} and \textit{F.~R. Olson}, Some theorems on \text{Bernoulli}
  and \text{Euler} numbers of higher order, Duke Mathematical Journal, vol.~21,
  no.~3, pp.~405--421  (1954).

\bibitem{cayley_00}
\textit{A.~Cayley}, On a theorem for the development of a factorial,
  Philosophical magazine, vol.~6, pp.~182--185  (1853).

\bibitem{cayley_01}
\textit{A.~Cayley}, On some numerical expansions, The Quarterly journal of pure
  and applied mathematics, vol.~3, pp.~366--369  (1860).

\bibitem{cayley_02}
\textit{A.~Cayley}, Note on a formula for \text{$\Delta^n 0^i/n^i$} when $n$,
  $i$ are very large numbers, Proceedings of the Royal Society of Edinburgh,
  vol.~14, pp.~149--153  (1887).

\bibitem{chabert_01}
\textit{J.-L. Chabert}, \textit{E.~Barbin}, \textit{M.~Guillemot},
  \textit{A.~Michel-Pajus}, \textit{J.~Borowczyk}, \textit{A.~Jebbar} and
  \textit{J.-C. Martzloff}, Histoire d'algorithmes: du caillou à la puce,
  Belin, France, 1994.

\bibitem{charalambides_01}
\textit{C.~A. Charalambides}, Enumerative Combinatorics, Chapman \& Hall/CRC,
  USA, 2002.

\bibitem{choi_01}
\textit{J.~Choi}, Certain integral representations of \text{Stieltjes}
  constants $\gamma_n$, Journal of Inequalities and Applications, 2013:532,
  pp.~1--10  (2013).

\bibitem{coffey_02}
\textit{M.~W. Coffey}, Series representations for the \text{Stieltjes}
  constants, Rocky Mountain Journal of Mathematics, vol.~44, pp.~443--477
  (2014), arXiv:0905.1111v2 (2009).
  
\bibitem{coffey_08}
\textit{M.~W. Coffey}, Addison-type series representation for the Stieltjes constants, 
Journal of Number Theory, vol.~130, pp.~2049--2064 (2010).

\bibitem{comtet_01}
\textit{L.~Comtet}, Advanced Combinatorics. The art of Finite and Infinite
  Expansions (revised and enlarged edition), D. Reidel Publishing Company,
  Dordrecht, Holland, 1974.

\bibitem{connon_06}
\textit{D.~F. Connon}, Some series and integrals involving the \text{Riemann}
  zeta function, binomial coefficients and the harmonic numbers. \text{Volume
  II(b)}, arXiv:0710.4024  (2007).

\bibitem{connon_01}
\textit{D.~F. Connon}, Some applications of the \text{Stieltjes} constants,
  arXiv:0901.2083  (2009).

\bibitem{connon_05}
\textit{D.~F. Connon}, A formula connecting the \text{Bernoulli} numbers with
  the \text{Stieltjes} constants, arXiv:1104.4772  (2011).

\bibitem{conway_01}
\textit{J.~H. Conway} and \textit{R.~K. Guy}, The Book of Numbers, Springer,
  New--York, 1996.

\bibitem{coppo_01}
\textit{M.-A. Coppo}, Nouvelles expressions des constantes de \text{Stieltjes},
  Expositiones Mathematic\ae, vol.~17, pp.~349--358  (1999).

\bibitem{copson_01}
\textit{E.~T. Copson}, Asymptotic Expansions, Cambridge University Press, Great
  Britain, 1965.

\bibitem{davis_02}
\textit{H.~T. Davis}, The approximation of logarithmic numbers, American
  Mathematical Monthly, vol.~64, no.~8, part II, pp.~11--18  (1957).
  
\bibitem{demoivre_01}
\textit{A.~\text{De Moivre}}, Miscellanea analytica de seriebus et quadraturis
  (with a supplement of 21 pages), J.~Thonson \& J.~Watts, Londini, 1730.

\bibitem{dilcher_01}
\textit{K.~Dilcher}, Generalized \text{Euler} constants for arithmetical
  progressions, Mathematics of Computation, vol.~59 , pp.~259--282  (1992).

\bibitem{dingle_01}
\textit{R.~B. Dingle}, Asymptotic Expansions: their Derivation and
  Interpretation, Academic Press, USA, 1973.

\bibitem{erdelyi_01}
\textit{A.~Erdélyi}, Asymptotic Expansions, Dover, USA, 1956.

\bibitem{ettingshausen_01}
\textit{A.~von Ettingshausen}, Die combinatorische Analysis als
  Vorbereitungslehre zum Studium der theoretischen höhern Mathematik,
  J.~B.~Wallishausser, Vienna, 1826.

\bibitem{euler_03}
\textit{L.~Euler}, Remarques sur un beau rapport entre les séries des
  puissances tant directes que réciproques, Histoire de l'Académie Royale des
  Sciences et Belles--Lettres, année MDCCLXI, Tome 17, pp.~83--106, A Berlin,
  chez Haude et Spener, Libraires de la Cour et de l'Académie Royale, 1768
  [read in 1749].

\bibitem{euler_02}
\textit{L.~Eulero}, Institutiones calculi differentialis cum eius usu in
  analysi finitorum ac doctrina serierum, Academiæ Imperialis Scientiarum
  Petropolitanæ, Saint--Petersburg, Russia, 1755.

\bibitem{ahmed_01}
\textit{L.~Fekih-Ahmed}, A new effective asymptotic formula for the
  \text{Stieltjes} constants, arXiv:1407.5567  (2014).

\bibitem{franel_01}
\textit{J.~Franel}, Note \no 245, L'Intermédiaire des mathématiciens, tome
  II, pp.~153--154  (1895).

\bibitem{furdui_01}
\textit{O.~Furdui}, \textit{M.~Bataille}, \textit{M.~Cibes},
  \textit{W.~Seaman}, \textit{D.~F. Connon}, \textit{K.~Lau} and
  \textit{J.~Lazzara}, Infinite sums and \text{Euler's} constant, The College
  Mathematics Journal, vol.~39, no.~1, pp.~71--72  (2008).

\bibitem{gauss_02}
\textit{C.~F. Gauss}, Disquisitiones generales circa seriem infinitam
  $1+\frac{\alpha\beta}{1\cdot\gamma}x+
  \frac{\alpha(\alpha+1)\beta(\beta+1)}{1\cdot2\cdot\gamma(\gamma+1)}xx+
  \frac{\alpha(\alpha+1)(\alpha+2)\beta(\beta+1)(\beta+2)}{1\cdot2\cdot3\cdot%
\gamma(\gamma+1)(\gamma+2)}x^3+\mathrm{etc}$, Commentationes Societatis Regiae
  Scientiarum Gottingensis recentiores, Classis Mathematic\ae, vol. II,
  pp.~3--46 [republished later in ``\text{Carl} \text{Friedrich} \text{Gauss}
  \text{Werke}'', vol.~3, pp.~265--327, \text{Königliche} \text{Gesellschaft}
  der \text{Wissenschaften}, \text{Göttingen}, 1866]  (1813).

\bibitem{gelfond_01}
\textit{A.~O. Gelfond}, The calculus of finite differences (3rd revised
  edition) [in Russian], Nauka, Moscow, USSR, 1967.

\bibitem{gerst_01}
\textit{Gerst}, Some series for \text{Euler's} constant, The American
  Mathematical Monthly, vol.~76, pp.~273--275  (1969).

\bibitem{gessel_01}
\textit{I.~Gessel} and \textit{R.~P. Stanley}, \text{Stirling} polynomials,
  Journal of Combinatorial Theory, vol.~A24, pp.~24--33  (1978).

\bibitem{glaisher_02}
\textit{G.~W.~L. Glaisher}, Congruences relating to the sums of products of the
  first $n$ numbers and to other sums of products, The Quarterly journal of
  pure and applied mathematics, vol.~31, pp.~1--35  (1900).

\bibitem{glaisher_01}
\textit{G.~W.~L. Glaisher}, On \text{Dr.~Vacca's} series for $\gamma$, The
  Quarterly journal of pure and applied mathematics, vol.~41, pp.~365--368
  (1910).

\bibitem{goldstine_01}
\textit{H.~H. Goldstine}, A History of Numerical Analysis from the 16th through
  the 19th Century, Springer--Verlag, New--York, Heidelberg, Berlin, 1977.

\bibitem{gould_01}
\textit{H.~W. Gould}, \text{Stirling} number representation problems,
  Proceedings of the American Mathematical Society, vol.~11, no.~3, pp.~447-451
   (1960).

\bibitem{gould_03}
\textit{H.~W. Gould}, An identity involving \text{Stirling} numbers, Annals of
  the Institute of Statistical Mathematics, vol.~17, \no.~1, pp.265--269
  (1965).

\bibitem{gould_02}
\textit{H.~W. Gould}, Note on recurrence relations for \text{Stirling} numbers,
  Publications de l'Institut Mathématique, Nouvelle série, vol.~6 (20),
  pp.~115--119  (1966).

\bibitem{knuth_01}
\textit{R.~L. Graham}, \textit{D.~E. Knuth} and \textit{O.~Patashnik}, Concrete
  mathematics: \text{A} foundation for computer science (2nd), Addison--Wesley,
  USA, 1994.

\bibitem{gram_01}
\textit{J.~P. Gram}, Note sur le calcul de la fonction $\zeta(s)$ de
  \text{Riemann}, Oversigt. K. Danske Vidensk. (Selskab Forhandlingar),
  pp.~305--308  (1895).

\bibitem{grunberg_01}
\textit{D.~B. Grünberg}, On asymptotics, \text{Stirling} numbers, gamma
  function and polylogs, Results in Mathematics, vol.~ 49, no.~1--2,
  pp.~89--125  (2006).
  
\bibitem{gunter_03_eng}
\textit{\text{N.~M.~Gunther (Günter)}} and \textit{\text{R.~O.~Kuzmin
  (Kusmin)}}, A Collection of Problems on Higher Mathematics. Vol.~3 (4th
  edition) [in Russian], Gosudarstvennoe izdatel'stvo tehniko--teoreticheskoj
  literatury, Leningrad, USSR, 1951.

\bibitem{hagen_01}
\textit{J.~G. Hagen}, Synopsis der höheren Analysis. Vol.~1. Arithmetische und
  algebraische Analyse, von Felix L.~Dames, Taubenstraße 47, Berlin, Germany,
  1891.

\bibitem{hardy_03}
\textit{G.~H. Hardy}, Note on \text{Dr.~Vacca's} series for $\gamma$, The
  Quarterly journal of pure and applied mathematics, vol.~43, pp.~215--216
  (1912).

\bibitem{hardy_02}
\textit{G.~H. Hardy}, Divergent series, Oxford at the Clarendan press, 1949.

\bibitem{hasse_01}
\textit{H.~Hasse}, Ein \text{Summierungsverfahren} für die \text{Riemannsche}
  $\zeta$-\text{Reihe}, Mathematische Zeitschrift, vol.~32, no.~1, pp.~458--464
   (1930).

\bibitem{hauss_01}
\textit{M.~Hauss}, Verallgemeinerte Stirling, Bernoulli und Euler Zahlen, deren
  Anwendungen und schnell konvergente Reihe für Zeta Funktionen
  (Ph.D.~dissertation), Aachen, Germany, 1995.

\bibitem{hayman_01}
\textit{W.~K. Hayman}, A generalisation of \text{Stirling's} formula, Journal
  für die reine und angewandte Mathematik, vol.~196, pp.~67--95  (1956).

\bibitem{hermite_01}
\textit{C.~Hermite}, Extrait de quelques lettres de \text{M.~Ch.~Hermite} à
  \text{M.~S.~Pincherle}, Annali di matematica pura ed applicata, serie III,
  tomo V, pp.~57--72  (1901).

\bibitem{hindenburg_01}
\textit{C.~F. Hindenburg}, Der polynomische Lehrsatz das wichtigste Theorem der
  ganzen Analysis nebst einigen Verwandten und andern Sätzen : Neu bearbeitet
  von \text{Tetens}, \text{Klügel}, \text{Kramp}, \text{Pfaff} und
  \text{Hindenburg}, bei Gerhard Fleischer, Leipzig, 1796.

\bibitem{howard_01}
\textit{F.~T. Howard}, Extensions of congruences of \text{Glaisher} and
  \text{Nielsen} concerning \text{Stirling} numbers, The Fibonacci Quarterly,
  vol.~28, no.~4, pp.~355--362  (1990).

\bibitem{howard_03}
\textit{F.~T. Howard}, Nörlund number \text{$B_n^{(n)}$}, in ``Applications of
  Fibonacci Numbers,'' vol.~5, pp.~355--366, Kluwer Academic, Dordrecht
  (1993).

\bibitem{howard_02}
\textit{F.~T. Howard}, Congruences and recurrences for \text{Bernoulli} numbers
  of higher orders, The Fibonacci Quarterly, vol.~32, no.~4, pp.~316--328
  (1994).

\bibitem{hwang_01}
\textit{H.~K. Hwang}, Asymptotic expansions for the \text{Stirling} numbers of
  the first kind, Journal of Combinatorial Theory, ser.~A 71, pp.~343--351
  (1995).

\bibitem{israilov_01}
\textit{M.~I. Israilov}, On the \text{Laurent} decomposition of
  \text{Riemann's} zeta function [in \text{Russian}], Trudy Mat. Inst. Akad.
  Nauk. SSSR, vol.~158, pp.~98--103  (1981).

\bibitem{jacobsthal_01}
\textit{E.~Jacobsthal}, Ueber die \text{Eulersche} konstante,
  Mathematisch--Naturwissenschaftliche Blätter, vol.~3, no.~9, pp.~153--154
  (1906).

\bibitem{jeffreys_02}
\textit{H.~Jeffreys} and \textit{B.~S. Jeffreys}, Methods of mathematical
  physics (second edition), University Press, Cambridge, Great Britain, 1950.

\bibitem{jensen_02}
\textit{J.~L. W.~V. Jensen}, Sur la fonction $\zeta(s)$ de \text{Riemann},
  Comptes-rendus hebdomadaires des séances de l'Académie des sciences, tome
  104, pp.~1156--1159  (1887).

\bibitem{jensen_04}
\textit{J.~L. W.~V. Jensen}, Opgaver til løsning \no 34, Nyt Tidsskrift for
  Matematik, Afdeling B, vol.~4, p.~54  (1893).

\bibitem{jensen_03}
\textit{J.~L. W.~V. Jensen}, Note \no 245. \text{Deuxième réponse}.
  \text{Remarques} relatives aux réponses du \text{MM. Franel et Kluyver},
  L'Intermédiaire des mathématiciens, tome II, pp.~346--347  (1895).

\bibitem{jordan_02}
\textit{C.~Jordan}, Sur des polynômes analogues aux polynômes de
  \text{Bernoulli}, et sur des formules de sommation analogues à celle de
  \text{MacLaurin--Euler}, Acta Scientiarum Mathematicarum (Szeged), vol.~4,
  no.~3-3, pp.~130--150  (1928--1929).

\bibitem{jordan_00}
\textit{C.~Jordan}, On \text{Stirling's Numbers}, Tohoku Mathematical Journal,
  First Series, vol.~37, pp.~254--278  (1933).

\bibitem{jordan_01}
\textit{C.~Jordan}, The calculus of finite differences, Chelsea Publishing
  Company, USA, 1947.

\bibitem{kenter_01}
\textit{F.~K. Kenter}, A matrix representation for \text{Euler's} constant
  $\gamma$, The American Mathematical Monthly, vol.~106, pp.~452--454  (1999).

\bibitem{kluyver_03}
\textit{J.~C. Kluyver}, On certain series of \text{Mr.~Hardy},
  Proc.~K.~Ned.~Akad. Wet., vol.~27, no.~3--4, pp.~314--323  (1924).

\bibitem{kluyver_02}
\textit{J.~C. Kluyver}, \text{Euler's} constant and natural numbers,
  Proc.~K.~Ned.~Akad. Wet., vol.~27, no.~1--2, pp.~142--144  (1924).

\bibitem{kluyver_01}
\textit{J.~C. Kluyver}, On certain series of \text{Mr.~Hardy}, The Quarterly
  journal of pure and applied mathematics, vol.~50, pp.~185--192  (1927).

\bibitem{knessl_02}
\textit{C.~Knessl} and \textit{M.~W. Coffey}, An asymptotic form for the
  \text{Stieltjes} constants $\gamma_k(a)$ and for a sum $s_\gamma(n)$
  appearing under the li criterion, Mathematics of Computation, vol.~80,
  no.~276, pp.~2197--2217  (2011).

\bibitem{knessl_01}
\textit{C.~Knessl} and \textit{M.~W. Coffey}, An effective asymptotic formula
  for the \text{Stieltjes} constants, Mathematics of Computation, vol.~80,
  no.~273, pp.~379--386  (2011).

\bibitem{knopp_01}
\textit{K.~Knopp}, Theory and applications of infinite series (2nd edition),
  Blackie \& Son Limited, London and Glasgow, UK, 1951.

\bibitem{knopp_02}
\textit{K.~Knopp}, Infinite sequences and series, Dover Publications Inc.,
  New-York, USA, 1956.

\bibitem{knuth_02}
\textit{D.~E. Knuth}, Two notes on notation, American Mathematical Monthly,
  vol.~99, no.~5, pp.~403--422  (1992).

\bibitem{korn_01}
\textit{G.~A. Korn} and \textit{T.~M. Korn}, Mathematical Handbook for
  Scientists and Engineers. Definitions, Theorems, and Formulas for Reference
  and Review (second, enlarged and revised edition), McGraw--Hill Book Company,
  New--York, 1968.

\bibitem{kowalenko_02}
\textit{V.~Kowalenko}, Generalizing the reciprocal logarithm numbers by
  adapting the partition method for a power series expansion, Acta
  Applicand\ae~Mathematic\ae, vol.~106, pp.~369--420  (2009).

\bibitem{kowalenko_01}
\textit{V.~Kowalenko}, Properties and applications of the reciprocal logarithm
  numbers, Acta Applicand\ae~Mathematic\ae, vol.~109, pp.~413--437  (2010).
  
\bibitem{skramer_01}
\textit{S.~Krämer}, Die Eulersche Konstante $\gamma$ und verwandte Zahlen
  (unpublished Ph.D.~manuscript, pers.~comm.), Göttingen, Germany, 2014.

\bibitem{kramp_01}
\textit{C.~Kramp}, Élements d'arithmétique universelle, L'imprimerie de
  Th.~F.~Thiriart, Cologne, 1808.

\bibitem{kratzer_01}
\textit{A.~Kratzer} and \textit{W.~Franz}, Transzendente Funktionen,
  Akademische Verlagsgesellschaft, Leipzig, Germany, 1960.

\bibitem{kreminski_01}
\textit{R.~Kreminski}, Newton--cotes integration for approximating
  \text{Stieltjes (generalized Euler)} constants, Mathematics of Computation,
  vol.~72, pp.~1379--1397  (2003).

\bibitem{krylov_01}
\textit{V.~I. Krylov}, Approximate calculation of integrals, The Macmillan
  Company, New-York, USA, 1962.

\bibitem{lagarias_01}
\textit{J.~C. Lagarias}, \text{Euler's} constant: \text{Euler's} work and
  modern developments, Bulletin (new series) of the American Mathematical
  Society, vol.~50 no.~4, pp.~527--628  (2013).

\bibitem{laplace_01}
\textit{P.-S. Laplace}, Traité de Mécanique Céleste. Tomes I--V, Courcier,
  Imprimeur--Libraire pour les Mathématiques, quai des Augustins, \no 71,
  Paris, France, 1800--1805.

\bibitem{laplace_02}
\textit{P.-S. Laplace}, Traité analytique des probabiblités, M\up{me} V\up{e}
  Courcier, Imprimeur--Libraire pour les Mathématiques, quai des Augustins,
  \no 57, Paris, France, 1812.

\bibitem{lavrik_01_eng}
\textit{A.~F. Lavrik}, On the main term of the divisor's problem and the power
  series of the \text{Riemann's} zeta function in a neighbourhood of its pole
  [in \text{Russian}], Trudy Mat. Inst. Akad. Nauk. SSSR, vol.~142,
  pp.~165--173  (1976).

\bibitem{lehmer_01}
\textit{D.~H. Lehmer}, \text{Euler} constants for arithmetical progressions,
  Acta Arithmetica, vol.~27, pp.~125--142  (1975).

\bibitem{lienard_01}
\textit{R.~Li\'enard}, Nombres de \text{Cauchy}, Intermédiaire des Recherches
  Mathématiques, vol.~2, no.~5, p.~38  (1946).

\bibitem{lindelof_01}
\textit{E.~Lindelöf}, Le calcul des résidus et ses applications à la
  théorie des fonctions, Gauthier--Villars, Imprimeur Libraire du Bureau des
  Longitudes, de l'École Polytechnique, Quai des Grands--Augustins, 55, Paris,
  1905.

\bibitem{louchard_01}
\textit{G.~Louchard}, Asymptotics of the \text{Stirling} numbers of the first
  kind revisited: \text{A} saddle point approach, Discrete Mathematics and
  Theoretical Computer Science, vol.~12, no.~2, pp.~167--184  (2010).

\bibitem{malgrange_01}
\textit{B.~Malgrange}, Sommation des séries divergentes, Expositiones
  Mathematic\ae, vol.~13, pp.~163--222  (1995).

\bibitem{mascheroni_01}
\textit{L.~Mascheronio}, Adnotationes ad calculum integralem \text{Euleri} in
  quibus nonnulla problemata ab \text{Eulero} proposita resolvuntur, Ex
  Typographia Petri Galeatii, Ticini, 1790.

\bibitem{mathstack_02}
\textit{Math.StackExchange}, A closed form for the series
  $\sum\frac{H^{2}_n - (\gamma+\ln n)^2}{n} $,
  http://math.stackexchange. \break com/questions/866382/  (2014).

\bibitem{matsuoka_01}
\textit{Y.~Matsuoka}, Generalized \text{Euler} constants associated with the
  \text{Riemann} zeta function, Number Theory and Combinatorics: Japan 1984,
  World Scientific, Singapore, pp.~279--295, 1985.

\bibitem{matsuoka_02}
\textit{Y.~Matsuoka}, On the power series coefficients of the \text{Riemann}
  zeta function, Tokyo Journal of Mathematics, vol.~12, no.~1, pp.~49--58
  (1989).

\bibitem{merlini_01}
\textit{D.~Merlini}, \textit{R.~Sprugnoli} and \textit{\text{M.~Cecilia}
  Verri}, The \text{Cauchy} numbers, Discrete Mathematics (Elsevier), vol.~306,
  pp.~1906--1920  (2006).

\bibitem{mezo_01}
\textit{I.~Mez\H{o}}, Gompertz constant, \text{Gregory} coefficients and a
  series of the logarithm function, Journal of Analysis \& Number Theory,
  vol.~2, no.~2, pp.~33--36  (2014).

\bibitem{mitrinovic_01}
\textit{D.~S. Mitrinovi\'c} and \textit{R.~S. Mitrinovi\'c}, Sur les nombres de
  \text{Stirling} et les nombres de \text{Bernoulli} d'ordre supérieur,
  Publications de la faculté d'électrotechnique de l'Université à
  Bélgrade, Série Mathématique et Physique, no. 43, pp.~1--63  (1960).

\bibitem{moser_01}
\textit{L.~Moser} and \textit{M.~Wyman}, Asymptotic development of the
  \text{Stirling} numbers of the first kind, Journal of the London Mathematical
  Society, vol.~s1--33, no.~2, pp.~133--146  (1958).

\bibitem{murray_01}
\textit{F.~J. Murray}, Formulas for factorial \text{$N$}, Mathematics of
  Computation, vol.~39, vol.~160, , pp.~655--662  (1982).

\bibitem{zhang_01}
\textit{Z.~Nan-You} and \textit{K.~S. Williams}, Some results on the
  generalized \text{Stieltjes} constant, Analysis, vol.~14, pp.~147--162
  (1994).

\bibitem{nemes_01}
\textit{G.~Nemes}, An asymptotic expansion for the \text{Bernoulli} numbers of
  the second kind, Journal of Integer Sequences, vol.~14, article 11.4.8
  (2011).

\bibitem{netto_01}
\textit{E.~Netto}, Lehrbuch der Combinatorik (2nd Edn.), Teubner, Leipzig,
  Germany, 1927.

\bibitem{nielsen_02}
\textit{N.~Nielsen}, En \text{Række} for \text{Eulers} \text{Konstant}, Nyt
  Tidsskrift for Matematik, Afdeling B, vol.~8, pp.~10--12  (1897).

\bibitem{nielsen_04}
\textit{N.~Nielsen}, Recherches sur les polynômes et les nombres de
  \text{Stirling}, Annali di Matematica Pura ed Applicata, vol.~10, no.~1,
  pp.~287--318  (1904).

\bibitem{nielsen_01}
\textit{N.~Nielsen}, Handbuch der \text{Theorie} der \text{Gammafunktion}, B.
  G. Teubner, Leipzig, Germany, 1906.

\bibitem{nielsen_03}
\textit{N.~Nielsen}, Recherches sur les polynômes de \text{Stirling},
  Hovedkommissionaer: Andr. Fred. Høst \& Søn, Kgl. Hof-Boghandel, Bianco
  Lunos Bogtrykkeri, København, Denmark, 1920.

\bibitem{norlund_02}
\textit{N.~E. Nörlund}, Vorlesungen über \text{Differenzenrechnung},
  Springer, Berlin, 1924.

\bibitem{norlund_01}
\textit{N.~E. Nörlund}, Sur les valeurs asymptotiques des nombres et des
  polynômes de \text{Bernoulli}, Rendiconti del Circolo Matematico di Palermo,
  vol.~10, no.~1, pp.~27--44  (1961).

\bibitem{olson_01}
\textit{F.~R. Olson}, Arithmetic properties of \text{Bernoulli} numbers of
  higher order, Duke Mathematical Journal, vol.~22, no.~4, pp.~641--653.
  (1955).

\bibitem{olver_01}
\textit{F.~W.~J. Olver}, Asymptotics and Special Functions, Academic Press,
  USA, 1974.

\bibitem{paplauskas_01_eng}
\textit{A.~B. Paplauskas}, Trigonometric series: from Euler to Lebesgue [in
  Russian], Nauka, Moscow, USSR, 1966.

\bibitem{pilehrood_01}
\textit{T.~H. Pilehrood} and \textit{K.~H. Pilehrood}, Criteria for
  irrationality of generalized \text{Euler's} constant, Journal of Number
  Theory, vol.~108, pp.~169--185  (2004).

\bibitem{polya_01_eng}
\textit{G.~P\'olya} and \textit{G.~Szeg\H{o}}, Problems and Theorems in
  Analysis I: Series, Integral calculus, Theory of functions, Springer--Verlag,
  Berlin, Germany, 1978.

\bibitem{proskuriyakov_01_eng}
\textit{I.~V. Proskuriyakov}, A Collection of Problems in Linear Algebra
  (fourth edition) [in Russian], Nauka, Moscow, USSR, 1970.

\bibitem{qi_01}
\textit{F.~Qi}, An integral representation, complete monotonicity, and
  inequalities of \text{Cauchy} numbers of the second kind, Journal of Number
  Theory, vol.~144, pp.~244--255  (2014).

\bibitem{ramis_01}
\textit{J.-P. Ramis}, Séries divergentes et procédés de resommation, École
  Polytechnique, Paris, France, 1991.

\bibitem{ramis_02}
\textit{J.-P. Ramis}, Séries divergentes et théories asymptotiques, Société
  mathématique de France, Paris, France, 1993.

\bibitem{rigaud_01}
\textit{S.~J. Rigaud}, Correspondence of scientific men of the seventeenth
  century, including letters of Barrow, Flamsteed, Wallis, and Newton, printed
  from the Originals [in 2 vols.], Oxford at the University Press, 1841.

\bibitem{riordan_01}
\textit{J.~Riordan}, An Introduction to Combinatorial Analysis, John Wiley \&
  Sons, Inc., USA, 1958.

\bibitem{rubinstein_01}
\textit{M.~O. Rubinstein}, Identities for the \text{Riemann} zeta function, The
  Ramanujan Journal, vol.~27, pp.~29--42  (2012).

\bibitem{rubinstein_02}
\textit{M.~O. Rubinstein}, Identities for the \text{Hurwitz} zeta function,
  \text{Gamma} function, and \text{$L$}-functions, The Ramanujan Journal,
  vol.~32, pp.~421--464  (2013).

\bibitem{eddin_01}
\textit{S.~Saad-Eddin}, Explicit upper bounds for the \text{Stieltjes}
  constants, Journal of Number Theory, vol.~133, no.~3, pp.~1027--1044  (2013).

\bibitem{eddin_02}
\textit{S.~Saad-Eddin}, On two problems concerning the
  \text{Laurent--Stieltjes} coefficients of \text{Dirichlet $L$--series} (Ph.D.
  thesis), University Lille 1, France, 2013.

\bibitem{salmieri_01}
\textit{A.~Salmeri}, Introduzione alla teoria dei coefficienti fattoriali, from
  ``Giornale di Matematiche di Battaglini'', vol.~90 (no.~10, serie 5),
  pp.~44--54  (1962).

\bibitem{sandham_01}
\textit{H.~F. Sandham}, Problem \no 4353: Euler's constant, American
  Mathematical Monthly, vol.~56, pp.~414  (1949).

\bibitem{sato_01}
\textit{H.~Sato}, On a relation between the \text{Riemann} zeta function and
  the \text{Stirling} numbers, Integers: Electronic Journal of Combinatorial
  Number Theory, vol.~8, no.~1  (2008).

\bibitem{schlaffli_01}
\textit{L.~Schläffli}, Sur les coëfficients du développement du produit
  $(1+x)(1+2x)\cdots \big(1+(n-1)\big)$ suivant les puissances ascendantes de
  $x$, Journal für die reine und angewandte Mathematik, vol.~43, pp.~1--22
  (1852).

\bibitem{schlaffli_02}
\textit{L.~Schläffli}, Ergänzung der abhandlung über die entwickelung des
  products $(1+x)(1+2x)\cdots \big(1+(n-1)\big)$ in band \text{XLIII} dieses
  journals, Journal für die reine und angewandte Mathematik, vol.~67,
  pp.~179--182  (1867).

\bibitem{schlomilch_04}
\textit{O.~Schlömilch}, Recherches sur les coefficients des facultés
  analytiques, Journal für die reine und angewandte Mathematik, vol.~44,
  pp.~344--355  (1852).

\bibitem{schlomilch_05}
\textit{O.~Schlömilch}, Compendium der höheren Analysis, Druck und Verlag von
  Friedrich Vieweg und Sohn, Braunschweig, Germany, 1853.

\bibitem{schlomilch_06}
\textit{O.~Schlömilch}, Compendium der höheren Analysis (2nd edn., in two
  volumes), Druck und Verlag von Friedrich Vieweg und Sohn, Braunschweig,
  Germany, 1861, 1866.

\bibitem{schlomilch_03}
\textit{O.~Schlömilch}, Nachschrift hierzu, Zeitschrift für angewandte
  Mathematik und Physik, vol.~25, pp.~117--119  (1880).

\bibitem{schroder_01}
\textit{E.~Schröder}, Bestimmung des infinitären \text{Werthes} des
  \text{Integrals} $\int\limits_0^1 (u)_n\, du$, Zeitschrift für angewandte
  Mathematik und Physik, vol.~25, pp.~106--117  (1880).

\bibitem{ser_01}
\textit{J.~Ser}, Sur une expression de la fonction $\zeta(s)$ de
  \text{Riemann}, Comptes-rendus hebdomadaires des séances de l'Académie des
  Sciences, Série 2, vol.~182, pp.~1075--1077  (1926).

\bibitem{shen_01}
\textit{L.-C. Shen}, Remarks on some integrals and series involving the
  \text{Stirling} numbers and $\zeta(n)$, The Transactions of the American
  Mathematical Society, vol.~347, no.~4, pp.~1391--1399  (1995).

\bibitem{shirai_01}
\textit{S.~Shirai} and \textit{K.~ichi Sato}, Some identities involving
  \text{Bernoulli} and \text{Stirling} numbers, Journal of Number Theory,
  vol.~90, pp.~130--142  (2001).

\bibitem{sondow_03}
\textit{J.~Sondow}, Analytic continuation of \text{Riemann's} zeta function and
  values at negative integers via \text{Euler's} transformation of series,
  Proceedings of the American Mathematical Society, vol.~120, no.~2,
  pp.~421--424  (1994).

\bibitem{srivastava_03}
\textit{H.~M. Srivastava} and \textit{J.~Choi}, Series Associated with the Zeta
  and Related Functions, Kluwer Academic Publishers, the Netherlands, 2001.

\bibitem{stamper_01}
\textit{P.~C. Stamper}, Table of \text{Gregory} coefficients, Mathematics of
  Computation, vol.~20, p.~465  (1966).
  
\bibitem{stankus_01_eng}
\textit{E.~Stankus}, A note about coefficients of \text{Laurent} series of the
  \text{Riemann} zeta function [in \text{Russian}], Zap.~Nauchn.~Sem.~LOMI,
  vol.~121, pp.~103--107  (1983).

\bibitem{stanley_01}
\textit{R.~P. Stanley}, Enumerative Combinatorics (1st Edn., 2nd printing),
  Cambridge University Press, 1997.

\bibitem{steffensen_01}
\textit{J.~F. Steffensen}, On \text{Laplace's} and \text{Gauss'}
  summation--formulas, Skandinavisk Aktuarietidskrift (Scandinavian Actuarial
  Journal), no.~1, pp.~1--15  (1924).

\bibitem{steffensen_02}
\textit{J.~F. Steffensen}, Interpolation (2nd Edn.), Chelsea Publishing
  Company, New--York, USA, 1950.

\bibitem{stirling_01}
\textit{J.~Stirling}, Methodus differentialis, sive Tractatus de summatione et
  interpolatione serierum infinitarum, Gul.~Bowyer, Londini, 1730.

\bibitem{tasaka_01}
\textit{T.~Tasaka}, Note on the generalized \text{Euler} constants,
  Mathematical Journal of Okayama University, vol.~36, pp.~29--34  (1994).

\bibitem{temme_02}
\textit{N.~M. Temme}, Asymptotic estimates of \text{Stirling} numbers, Studies
  in Applied Mathematics, vol.~89, pp.~233--243  (1993).

\bibitem{timashev_01}
\textit{A.~N. Timashev}, On asymptotic expansions of \text{Stirling} numbers of
  the first and second kinds, Discrete Mathematics and Applications, vol.~8,
  no.~5, pp.~533--544  (1998).
  
\bibitem{todd_01}
\textit{J.~J.~Y. Liang} and \textit{J.~Todd}, The \text{Stieltjes} constants,
  Journal of Research of the National Bureau of Standards---Mathematical
  Sciences, vol.~76B, nos.~3--4, pp.~161--178  (1972).

\bibitem{tricomi_01}
\textit{F.~G. Tricomi} and \textit{A.~Erdélyi}, The asymptotic expansion of a
  ratio of gamma functions, Pacific Journal of Mathematics, vol.~1, no.~1,
  pp.~133--142  (1951).

\bibitem{newton_01}
\textit{H.~W. Turnbull}, The correspondence of Isaac Newton [vols.~1--7], Royal
  Society at the University Press, Cambridge, 1959--1977.

\bibitem{tweedie_01}
\textit{C.~Tweedie}, The \text{Stirling} numbers and polynomials, Proceedings
  of the Edinburgh Mathematical Society, vol.~37, pp.~2--25  (1918).

\bibitem{vacca_01}
\textit{G.~Vacca}, A new series for the \text{Eulerian} constant, The Quarterly
  journal of pure and applied mathematics, vol.~41, pp.~363--364  (1910).
  
\bibitem{van_veen_01}
\textit{S.~C. \text{Van Veen}}, Asymptotic expansion of the generalized
  \text{Bernoulli} numbers \text{$B_n^{(n-1)}$} for large values of $n$ ($n$
  integer), Indagationes Mathematic\ae~(Proceedings of the Koninklijke
  Nederlandse Akademie van Wetenschappen. Series A, Mathematical sciences),
  vol.~13, pp.~335--341  (1951).

\bibitem{vorobiev_01}
\textit{N.~N. Vorobiev}, Theory of series (4th edition, enlarged and revised)
  [in Russian], Nauka, Moscow, USSR, 1979.

\bibitem{wachs_01}
\textit{S.~Wachs}, Sur une propriété arithmétique des nombres de
  \text{Cauchy}, Bulletin des Sciences Mathématiques, deuxi\`eme s\'erie.
  vol.~71, pp.~219--232  (1947).

\bibitem{watson_01}
\textit{G.~N. Watson}, An expansion related to \text{Stirling's} formula,
  derived by the method of steepest descents, The Quarterly journal of pure and
  applied mathematics, vol.~48, pp.~1--18  (1920).

\bibitem{weisstein_04}
\textit{E.~W. Weisstein}, CRC Concise Encyclopedia of Mathematics (2nd Edn.),
  Chapman \& Hall/CRC, USA, 2003.
  
\bibitem{weisstein_06}
\textit{E.~W. Weisstein}, Riemann zeta function zeros, http://mathworld.wolfram.com/RiemannZeta \break
 FunctionZeros.html, 2015.

\bibitem{whittaker_01}
\textit{E.~Whittaker} and \textit{G.~N. Watson}, A course of modern analysis.
  An introduction to the general theory of infinite processes and of analytic
  functions, with an account of the principal transcendental functions (third
  edition), Cambridge at the University Press, Great Britain, 1920.

\bibitem{wilf_02}
\textit{H.~S. Wilf}, The asymptotic behavior of the \text{Stirling} numbers of
  the first kind, Journal of Combinatorial Theory, ser.~A 64, pp.~344--349
  (1993).

\bibitem{wilf_01}
\textit{H.~S. Wilf}, Generatingfunctionology (2nd), Academic Press, Inc., USA,
  1994.

\bibitem{xia_01}
\textit{L.~Xia}, The \text{parameterized-Euler-constant} function
  $\gamma_a(z)$, Journal of Number Theory, vol.~133, no.~1, pp.~1--11  (2013).

\bibitem{young_01}
\textit{P.~T. Young}, A 2-adic formula for \text{Bernoulli} numbers of the
  second kind and for the \text{Nörlund} numbers, Journal of Number Theory,
  vol.~128, pp.~2951--2962  (2008).

\bibitem{zhao_01}
\textit{F.-Z. Zhao}, Sums of products of \text{Cauchy} numbers, Discrete
  Mathematics, vol.~309, pp.~3830--3842  (2009).

\end{thebibliography}
\end{document}